\documentclass[10pt]{article}
\usepackage{amsmath,amsfonts,amstext,amsthm,epsfig}
\begin{document}
\font\tenrm=cmr10

\newtheorem{Theorem}{Theorem}
\newtheorem{Lemma}{Lemma}
\newtheorem{Corol}{Corollary}

\newtheorem{Prop}{Proposition}
\newtheorem{Rem}{Remark}
\makeatletter \addtocounter{section}{0}
\renewcommand{\theequation}{\thesection.\arabic{equation}}
\@addtoreset{equation}{section} \makeatother
\def\R{\mathbb R}
\def\N{\mathbb N}

\title{\bf Continuous Dependence for Backward Parabolic Operators with
Log-Lipschitz Coefficients}

\author{ Daniele Del Santo {\small and}  Martino Prizzi\\ {\small Dipartimento di
Matematica e Informatica, Universit\`a di Trieste}\\[-0.1 cm]
{\small Via A.~Valerio 12/1, 34127 Trieste, Italy} }

\date{\today}

\maketitle

\begin{abstract}

We prove continuous dependence on Cauchy data for a backward
parabolic operator whose coefficients are Log-Lipschitz continuous
in time.\vskip 0.3 cm\noindent {\bf Keywords}: backward parabolic
operator, ill posed problem, continuous dependence on data. \vskip
0.3 cm\noindent  {\bf 2000 MSC}: 35K10, 34A55, 35B30, 35R25.

\end{abstract}

\section{Introduction}

In this paper we prove a new continuous dependence result for
solutions of the Cauchy problem associated to the backward
parabolic operator
\begin{equation}
P=\partial_t
+\sum_{i,j}\partial_{x_i}(a_{i,j}(t,x)\partial_{x_j})+\sum_j
b_j(t,x)\partial_{x_j}+c(t,x) \label{intro1}
\end{equation}
on the strip $ [0,T]\times {\mathbb R}^n$.

It is well known that the Cauchy problem for (\ref{intro1}), when the data are given on $\{t=0\}$ and the matrix $(a_{i,j})^n_{i,j=1}$ is supposed to be symmetric and positive definite, is an
ill-posed problem: due to the smoothing effect of forward parabolic
operators, the existence of the solutions is 
not ensured for all choice of data. 
Concerning uniqueness, we can say that an important role is played by the functional space in which the uniqueness property is looked for. In fact a classical result of  Tychonoff  in \cite{Ty} proves that there exists
 a function $u\in C^\infty (\R\times\R^n)$
satisfying $\partial_t u-\Delta u\equiv 0$ in $\R\times\R^n$,
$u(0,\cdot) \equiv 0$ in $\R^n$, but $u\not\equiv 0$ in all open
subset of $\R\times\R^n$. On the other hand, in
\cite{LM}, Lions and Malgrange proved that $P$ enjoys the
uniqueness property in ${\mathcal H}_1:=H^1([0,T], L^2(\R^n))\cap
L^2([0,T], H^{2}(\R^n))$, provided the coefficients $a_{i,j}$'s are
sufficiently smooth with respect to $x$ and Lipschitz continuous
with respect to $t$.

If one considers the  Cauchy problem for (\ref{intro1}) as an
inverse problem (namely: a final time problem) for a forward
parabolic operator (see \cite[Ch. 3]{I}), it turns out that uniqueness is a very weak
property. Indeed it furnishes only a qualitative feature of the
solutions and gives no useful information for computational
purposes. 

In his celebrated paper \cite{J}, John introduced the
notion of {\it well-behaved problem}, which is now typical in the
context of ill-posed problems. According to John a problem is {\it
well-behaved} if
{\it ``only a fixed percentage of the significant
digits need be lost in determining the solution from the data"} \cite[p. 552]{J}. More precisely we may say that  a problem is 
well-behaved if its solutions in a space $\mathcal H$ depend
H\"older continuously on the data belonging to a space $\mathcal
K$, provided they satisfy a prescribed bound. 

In their paper \cite{AN}, Agmon and Nirenberg proved, among other
things, that the Cauchy problem for (\ref{intro1}) is well-behaved
in ${\mathcal E}:=C^0([0,T], L^2({\mathbb R}^n))\cap C^0([0,T[,
H^1({\mathbb R}^n))\cap C^1([0,T[, L^2({\mathbb R}^n))$ with data
in $L^2(\R^n)$, provided the coefficients $a_{i,j}$'s are
sufficiently smooth with respect to $x$ and Lipschitz continuous
with respect to $t$. In order to achieve their result, which is
stated in a very general and  abstract setting, they developed the so
called {\it logarithmic convexity technique}. The main step
consists in proving that the function $t\mapsto
\log\|u(t,\cdot)\|_{L^2}$ is convex for every solution
$u\in{\mathcal E }$ of (\ref{intro1}). In the same year Glagoleva
\cite{G} obtained essentially the same result for a concrete
operator like (\ref{intro1}) with time independent coefficients.
Her proof rests on energy estimates obtained through integration
by parts. Some years later Hurd \cite{H} developed the technique
of Glagoleva so as to cover the case of a general operator of type
(\ref{intro1}), with coefficients depending Lipschitz continuously
on time. The results of \cite{AN,G,H} can be summarized as
follows:

\bigskip

{\it For every $T'\in\,\,]0,T[$ and $D>0$ there exist $\rho>0$,
$0<\delta<1$ and $M>0$ such that, if $u\in {\mathcal E}$ is a
solution of $Pu\equiv0$ on $[0,T]$ with
$\|u(0,\cdot)\|_{L^2}\leq\rho$ and $\|u(t,\cdot)\|_{L^2}\leq D$ on
$[0,T]$, then $$\sup_{t\in[0,T']}\|u(t,\cdot)\|_{L^2}\leq
M\|u(0,\cdot)\|_{L^2}^\delta.$$  The constants $\rho$, $M$ and
$\delta$ depend only on $T'$ and $D$, on the ellipticity constant
of $P$, on the $L^\infty$ norms of the coefficients $a_{i,j}$'s,
$b_i$'s, $c$ and of their spatial derivatives, and on the Lipschitz
constant of the coefficients $a_{i,j}$'s with respect to time.}

\bigskip

In \cite{LM,AN,H}, Lipschitz continuity of the coefficients
$a_{i,j}$'s with respect to time plays an essential role. The
possibility of replacing Lipschitz continuity by simple continuity
was ruled out by Miller \cite{Mil} and more recently by Mandache
\cite{M}. They constructed examples of operators of the form
(\ref{intro1}) which do not enjoy the uniqueness property in
${\mathcal H}_1$. In the example of Miller the coefficients
$a_{i,j}$'s are H\"older continuous in time, while in the more
refined example of Mandache the modulus of continuity $\bar\mu$ of
the coefficients $a_{i,j}$'s with respect to time is such that
$\int_0^1(1/\bar\mu(s))ds<+\infty$. On the other hand, in
\cite{DSP} the authors of the present paper proved that, if
$\bar\mu$ satisfies the {\it Osgood condition }
$\int_0^1(1/\bar\mu(s))ds=+\infty$, then the operator $P$ enjoys
the uniqueness property in ${\mathcal H}_1$. Therefore it would be
natural to conjecture that if the Osgood condition  is satisfied,
then the Cauchy problem for (\ref{intro1})  is well-behaved in
$\mathcal E$ with data in $L^2(\R^n)$. Unfortunately this is not
true. Let $\mu(s):=s(1+|\log s|)$. A function whose modulus of
continuity is $\mu$ is called {\it Log-Lipschitz continuous}. Obviously $\mu$ satisfies the Osgood condition. In
the Appendix, we show that it is possible to construct:

\begin{itemize}

\item a sequence $(L_n)_{n\in\N}$ of backward uniformly parabolic
operators with space-periodic uniformly Log-Lipschitz continuous
coefficients in the principal part and space-periodic uniformly
bounded coefficients in lower order terms;

\item a sequence $(u_n)_{n\in\N}$ of space-periodic smooth
uniformly bounded solutions of $L_nu_n=0$ on $[0,1]\times\R^2$;

\item a sequence $(t_n)_{n\in\N}$ of real numbers, with $t_n\to
0$ as $n\to\infty$;

\end{itemize}
such that
$$\lim_{n\to\infty}\|u_n(0,\cdot)\|_{L^2([0,2\pi]\times[0,2\pi])}=0
$$ and $$ \lim_{n\to\infty} {\|u_n(t_n, \cdot,
\cdot)\|_{L^2([0,2\pi] \times [0,2\pi])}\over \|u_n(0, \cdot,
\cdot)\|^\delta_{L^2([0,2\pi] \times [0,2\pi])}}=+\infty, $$ for
every $\delta>0$.

Therefore it is not possible to obtain a result similar to that of
Hurd  or Agmon and Nirenberg  if Lipschitz continuity is replaced
by Log-Lipschitz continuity. 

If the coefficients
$a_{i,j}$'s are Log-Lipschitz continuous in time, we are able to
prove a weaker continuous dependence result. Our main result can
be stated as follows:

\bigskip

{\it For every $T'\in\,\,]0,T[$ and $D>0$ there exist $\rho>0$,
$0<\delta<1$ and $M,N>0$ such that, if $u\in {\mathcal E}$ is a
solution of $Pu\equiv0$ on $[0,T]$ with
$\|u(0,\cdot)\|_{L^2}\leq\rho$ and $\|u(t,\cdot)\|_{L^2}\leq D$ on
$[0,T]$, then $$\sup_{t\in[0,T']}\|u(t,\cdot)\|_{L^2}\leq M
e^{-N|\log \|u(0,\cdot)\|_{L^2}|^{\delta }}.$$  The constants
$\rho$, $M$, $N$ and $\delta$ depend only on $T'$ and $D$, on the
ellipticity constant of $P$, on the $L^\infty$ norms of the
coefficients $a_{i,j}$'s, $b_i$'s, $c$ and of their spatial
derivatives, and on the Log-Lipschitz constant of the coefficients
$a_{ij}$'s with respect to time.}

\bigskip
\noindent
As a consequence, going back to John's terminology, if one denotes by $\phi(n)$ the number of digits of the $L^2$ norm
of the data which are necessary to determine $n$ digits of the $L^2$ norm of the solution, one has that $\phi(n)$ grows at most polynomially in $n$.

\smallskip

Our proof relies on weighted energy estimates similar to
those of Glagoleva and Hurd. In order to overcome the obstructions
created by the lack of time differentiability of the coefficients
$a_{i,j}$'s, we exploit a microlocal approximation procedure
originally developed by Colombini and Lerner in \cite{CL} in the study of the Cauchy problem for hyperbolic operators having Log-Lipschitz coefficients.

The plan of the paper is the following. In Section 2 we introduce
notations and we state our results: Theorem 1 contains the
weighted energy estimates that we mentioned above; Theorem 2 is a
local continuous dependence result; Theorem 3 is a global
continuous dependence result. Section 3 is devoted to the proof of
Theorem 1, while Section 4 is devoted to the proofs of Theorems 2
and 3. Finally, in the Appendix we outline the construction of a
counterexample to H\"older continuous dependence.

 \section{Results}
\subsection{Notations}
We consider the following backward parabolic equation
\begin{equation}
\partial_t u +\sum_{i,j}\partial_{x_i}(a_{i,j}(t,x)\partial_{x_j}u)+\sum_j b_j(t,x)\partial_{x_j}u+c(t,x)u=0
\label{1}
\end{equation}
on the strip $ [0,T]\times {\mathbb R}^n \ni (t,x)$. We suppose that
\begin{itemize}
\item for all $(t,x)\in [0,T]\times {\mathbb R}^n$ and for all $i,j=1\dots n$, 
$$
a_{i,j}(t,x)=a_{j,i}(t,x);
$$
\item there exists $k>0$ such that, for all $(t,x, \xi)\in [0,T]\times {\mathbb R}^n\times {\mathbb R}^n$,
$$ 
k|\xi|^2\leq \sum_{i,j}a_{i,j}(t,x)\xi_i\xi_j\leq k^{-1}|\xi|^2;
$$
\item for all $i,j=1,\dots , n$, $a_{i,j}\in {\rm LogLip}([0,T], L^\infty({\mathbb R}^n))\cap L^\infty([0,T], C^2_b({\mathbb R}^n))$ and $b_j$, $c\in L^\infty([0,T], C^2_b({\mathbb R}^n))$.
\end{itemize}
We set
$$
\begin{array}{ll}
&\displaystyle{A_{LL}:= \sup\,\{{|a_{i,j}(t,x)-a_{i,j}(s,x)|\over |t-s|(1+|\log|t-s||)}
\,:\,i,j=1,\dots,n,\  x\in {\mathbb R}^n,}\\[0.3 cm] 
&\qquad\qquad\qquad\qquad\qquad\qquad\qquad\qquad\qquad\quad \displaystyle{t,s\in[0,T],\ 0<|t-s|\leq 1\},} \\[0.3 cm] 
&A:=\sup\,\{|\partial_x^\alpha a_{i,j}(t,x)|\,:\,i,j=1,\dots,n,\ \alpha\in {\mathbb N}^n, \ |\alpha|\leq 2,\\[0.3 cm] 
&\qquad\qquad\qquad\qquad\qquad\qquad\qquad\qquad\qquad\qquad\qquad\quad
(t,x)\in[0,T]\times {\mathbb R}^n\},\\[0.3 cm]
&B:=\sup\,\{|\partial_x^\alpha b_{j}(t,x)|\,:\,j=1,\dots,n,\ \alpha\in {\mathbb N}^n, \ |\alpha|\leq 2,\\[0.3 cm] 
&\qquad\qquad\qquad\qquad\qquad\qquad\qquad\qquad\qquad\qquad\qquad\quad
(t,x)\in[0,T]\times {\mathbb R}^n\},\\[0.3 cm]
&C:=\sup\,\{|\partial_x^\alpha c(t,x)|\,:\,\alpha\in {\mathbb N}^n, \ |\alpha|\leq 2, \ 
(t,x)\in[0,T]\times {\mathbb R}^n\}.
\end{array}
$$

\subsection{Weight function}

For $s>0$, let  $\mu(s)=s(1+|\log s|)$. For $\tau \geq 1$ we define 
$$
\theta(\tau):=\int_{1/\tau}^1{1\over \mu(s)}\, ds=\log(1+\log\tau).
$$
The function $\theta:[1,+\infty[\,\to[0,+\infty[$ is bijective and strictly increasing.
For $y\in\, ]0,1]$ and $\lambda>1$, we set $\psi_\lambda(y):=\theta^{-1}(-\lambda\log y)=\exp(y^{-\lambda}-1)$ and we define
$$
\Phi_\lambda(y):=\int_1^y\psi_\lambda(z)\, dz.
$$
The function $\Phi_\lambda:\,]0,1]\to\,]-\infty, 0]$ is bijective and strictly increasing; moreover it is easy to verify that it satisfies 
\begin{equation}
y\, \Phi_\lambda''(y) =-\lambda(\Phi_\lambda'(y))^2\mu({1\over \Phi_\lambda'(y)})=-\lambda\, \Phi_\lambda'(y) (1+|\log {1\over \Phi_\lambda'(y)}|).
\label{ode}
\end{equation}
We collect in the following lemma, the proof of which is left to the reader, some interesting and elementary properties
of the functions $\psi_\lambda $ and $\Phi_\lambda$.
\begin{Lemma}
\label{lemma1}
Let $\zeta>1$. Then, for $y\leq 1/\zeta$, 
\begin{equation}
\psi_\lambda(\zeta y)=\exp{(\zeta^{-\lambda}-1)}(\psi_\lambda(y))^{\zeta^{-\lambda}}.
\label{a}
\end{equation}

Define $\Lambda_\lambda(y):=y\,\Phi_\lambda(1/y)$. Then the function $\Lambda_\lambda:[1,+\infty[\,\to\,]-\infty, 0]$ is bijective  and 
\begin{equation}
\lim_{z\to -\infty} -{1\over z}\,\psi_\lambda({1\over \Lambda_\lambda^{-1}(z)})=+\infty.
\label{b}
\end{equation}
\end{Lemma}

\subsection{Main results}

Let ${\mathcal E}:=C^0([0,T], L^2({\mathbb R}^n))\cap C^0([0,T[, H^1({\mathbb R}^n))\cap C^1([0,T[, L^2({\mathbb R}^n))$.

\begin{Theorem}
\label{theorem1}
There exist $\bar \lambda>1$, $\alpha_1$, $\bar \gamma$, $M>0$ depending only on $A_{LL}$, $A$, $B$, $C$, 
$k$, $T$ such that, setting $\alpha:= \max\{\alpha_1, 1/T\}$, $\sigma:=1/\alpha$ and choosing $\tau\in \,]0,\sigma/2[$, if $\beta\geq \sigma+\tau$, $\lambda\geq\bar \lambda$, $\gamma\geq\bar \gamma$ and if $u\in {\mathcal E}$ is a solution of the equation (\ref{1}), then 
\begin{equation}
\label{2}
\begin{array}{l}
\displaystyle{\int_0^s e^{2\gamma t} e^{-2\beta\Phi_\lambda((t+\tau)/ \beta)} 
\|u(t,\cdot)\|^2_{H^{1-\alpha t}}\, dt}\\[0.3 cm]
\qquad\qquad\qquad\displaystyle{ \leq M \Big((s+\tau)e^{2\gamma s} e^{-2\beta\Phi_\lambda((s+\tau)/ \beta)} 
\|u(s,\cdot)\|^2_{H^{1-\alpha s}}}\\[0.3 cm]
\qquad\qquad\qquad\qquad\qquad\qquad\qquad
\displaystyle{+\tau\, \Phi'_\lambda({\tau/ \beta})e^{-2\beta\Phi_\lambda({\tau/\beta})}
\|u(0,\cdot)\|^2_{L^2}\Big),}
\end{array}\end{equation}
for all $0\leq s\leq\sigma$.
\end{Theorem}

\begin{Theorem}
\label{theorem2}
There exist $\sigma>0$ such that for all $\, \bar \sigma\in \, ]0,\sigma/4[\,$ there exist  $ \rho$, $\bar M$, $ N$, $ \delta>0$ such that, if  $u\in {\mathcal E}$ is a solution of the equation (\ref{1}) with $\|u(0,\cdot)\|_{L^2}\leq \rho$, then
\begin{equation}
\sup_{t\in [0,\bar\sigma]}\|u(t,\cdot)\|_{L^2}\leq \bar M(1+\|u(\sigma,\cdot)\|_{L^2})e^{- N(|\log \|u(0,\cdot)\|_{L^2}|)^{\delta}}.
\label{3}
\end{equation}
\end{Theorem}

\begin{Theorem}
\label{theorem3}
For all $\, T'\in \,]0,T[$  and for all $D>0$ there exist $ \rho'$, $ M'$, $ N'$, $ \delta'>0$, depending only on 
$A_{LL}$, $A$, $B$, $C$,  $k$, $T$, $T'$, $D$, such that if $u\in {\mathcal E}$ is a solution of the equation (\ref{1}) with $\sup_{t\in [0,T]}\|u(t,\cdot)\|_{L^2}\leq D$ and $\|u(0,\cdot)\|_{L^2}\leq \rho'$, then
\begin{equation}
\sup_{t\in [0,T']}\|u(t,\cdot)\|_{L^2}\leq M' e^{-N'|\log \|u(0,\cdot)\|_{L^2}|^{\delta '}}.
\label{4}
\end{equation}
\end{Theorem}

\section{Proof of Theorem 1}
\subsection{Dyadic decomposition}
We collect here some well known facts on the Littlewood-Paley dyadic decomposition, referring to \cite{Bo} and \cite{CL} for the details. Let $\varphi\in C^\infty({\mathbb R})$, $\varphi(x)=1$ if $x\leq 1$, $\varphi(x)=0$ if $x\geq 2$,  $\varphi$   decreasing. We set $\varphi_0(\xi)=\varphi(|\xi|)$ and,  if $\nu$ is an integer greater than or equal to 1, $\varphi_\nu(\xi)=\varphi_0(\xi/2^\nu)-\varphi_0(\xi/2^{\nu-1})$. Let $w$ be a tempered distribution in ${H}^{-\infty}({\mathbb R}^n)$; we define
\[
\begin{array}{l}
\displaystyle{w_\nu(x)=\varphi_\nu(D_x)w(x)= {1\over (2\pi)^n}\int e^{ix \cdot \xi}\varphi_\nu(\xi)\hat w(\xi)\, d\xi}\\[0.3 cm]
\qquad\qquad\qquad\qquad\qquad\qquad\qquad\qquad\displaystyle{={1\over (2\pi)^n}\int\hat{\varphi_\nu}( y)w(x-y)\, dy}.
\end{array}
\]
For all $\nu$, $w_\nu$ is an entire analytic function belonging to ${L}^2$. We have
\begin{itemize}
\item for all $\nu\geq 1$
\begin{equation}
\label{d1}
2^{\nu-1}\|w_\nu\|_{L^2}\leq \|\nabla_x w_\nu\|_{(L^2)^n}\leq 2^{\nu+1}  \| w_\nu\|_{L^2},
\end{equation}
where the inequality on the right hand side holds also for $\nu=0$;
\item there exist $K$ such that, for all $s\in [0,1]$, 
\begin{equation}
\label{d2}
K \sum_{\nu=0}^{+\infty}2^{2s\nu}\|w_\nu\|^2_{ L^2} \leq \|w\|^2_{{ H}^s}\leq {1\over K} \sum_{\nu=0}^{+\infty}2^{2s\nu}\|w_\nu\|^2_{{ L}^2};
\end{equation}
\item if the function $u:[0,T[\,\to L^2({\mathbb R}^n)$ is of class $C^1$, then the function
$u_\nu: [0,T[\,\to C^m_b({\mathbb R}_x^n)\cap H^s({\mathbb R}_x^n)$ is of class $C^1$ for all $s\geq 0$ and for all $m\in {\mathbb N}$ and,  for all $s\geq 0$ and for all $\alpha\in{\mathbb N}^n$,
$$
\partial_t\partial_x^\alpha u_\nu= \partial_x^\alpha\partial_t u_\nu\in 
C^0([0,T[\times {\mathbb R}^n)\cap C^0([0,T[, L^2({\mathbb R}^n));
$$
\item
if $a\in C^2_b({\mathbb R}^n)$, then there exits $Q>0$ such that, for all $\nu, \mu\in {\mathbb N}$,
$$
\|[\varphi_\nu(D_x), a]\varphi_\mu(D_x)\|_{{\mathcal L}(L^2, L^2)}
\leq \left\{
\begin{array}{ll}
\displaystyle{Q\, 2^{-2\nu}}&\quad{\rm if}\quad |\mu-\nu|\leq 2,\\[0.5 cm]
\displaystyle{Q\, 2^{-2\max\{\nu, \mu\}}}\ &\quad{\rm if}\quad |\mu-\nu|\geq 3,
 \end{array}
 \right.
 $$
 where $[\varphi_\nu(D_x), a]w(x)=(\varphi_\nu(D_x)(aw))(x)-a(x)(\varphi_\nu(D_x)u)(x)$ is the commutator between $\varphi _\nu(D_x)$ and $a$, 
 $\|\cdot\|_{{\mathcal L}(L^2, L^2)}$ denotes the norm operator in $L^2$ to $L^2$ and the constant $ Q$ depends only on $\|a\|_{C^2_b}$.
\end{itemize}

\subsection{Preliminaries}

Let $u(t,x)\in  {\mathcal E}$ be a solution of the equation (\ref{1}).
We set 
$$
w(t,x):= e^{\gamma t} e^{-\beta \Phi_\lambda ({t+\tau\over \beta})}u(t,x),
$$
$u_\nu(t,x):=\varphi_\nu(D_x)u(t,x)$, $w_\nu(t,x):=\varphi_\nu(D_x)w(t,x)$,
$v_\nu(t,x):= 2^{-\alpha\nu t}w_\nu(t,x)$, where the constants $\alpha$, $\lambda$ and $\gamma$ will be determined later, $\sigma:=1/\alpha$,  $\tau$ is chosen in $]0,\sigma/2[$ and $ \beta\geq \sigma+\tau$.
The function $v_\nu$ satisfies
\begin{equation}
\label{eq2}
\begin{array}{ll}
\partial_t v_\nu=\displaystyle{\gamma v_\nu-\sum_{i,j}\partial_{x_i}(a_{i,j}(t,x)\partial_{x_j} v_\nu)-\Phi'_\lambda({t+\tau\over \beta})v_\nu-\alpha(\log 2)\nu v_\nu}\\[0.3cm]
\ \ \quad\qquad\qquad\qquad\quad\qquad\qquad\displaystyle{-\sum_j b_j(t,x)\partial_{x_j}v_\nu-c(t,x)v_\nu+X_\nu(t,x),}
\end{array}
\end{equation}
where
$$
\begin{array}{ll}
X_\nu(t,x):=\displaystyle{ -\sum_{i,j} \partial_{x_i}([\varphi_\nu(D_x),a_{i,j}(t,x)]2^{-\alpha \nu t}\partial_{x_j} w)}\\[0.3 cm]
\qquad\qquad\qquad\quad\displaystyle
{-\sum_{j} [\varphi_\nu(D_x),b_j(t,x)]2^{-\alpha \nu t}\partial_{x_j} w-[\varphi_\nu(D_x),c(t,x)]2^{-\alpha \nu t}w.}
\end{array}
$$
Setting $A(t,x):= (a_{i,j}(t,x))^n_{i,j=1}\in {\mathcal M}^{n\times n}$ and $B(t,x):=
(b_{j}(t,x))^n_{j=1}\in {\mathbb R}^n$, we compute the scalar product of (\ref{eq2})
with $(t+\tau)\partial_t v_\nu$ and we obtain
\begin{equation}
\label{****}
\begin{array}{ll}
\displaystyle{(t+\tau)\|\partial_t v_\nu(t, \cdot)\|^2_{L^2}=}&\displaystyle{\gamma(t+\tau)\langle v_\nu(t,\cdot), \partial_tv_\nu(t,\cdot)\rangle_{L^2}}\\[0.3 cm]
&\displaystyle{\ +(t+\tau)\langle A(t,\cdot)\nabla_x v_\nu(t,\cdot), \nabla_x\partial_t v_\nu(t,\cdot)\rangle_{(L^2)^n}}\\[0.3 cm]
&\displaystyle{\ \  -(t+\tau)\Phi_\lambda'({t+\tau\over \beta})\langle v_\nu(t,\cdot), \partial_tv_\nu(t,\cdot)\rangle_{L^2}}\\[0.3 cm]
&\displaystyle{\ \ \ -\alpha(\log 2)\nu(t+\tau)\langle v_\nu(t,\cdot), \partial_t v_\nu (t,\cdot)\rangle_{L^2}}\\[0.3 cm]
&\displaystyle{\ \ \ \ -(t+\tau)\langle B(t,\cdot)\cdot \nabla_xv_\nu(t,\cdot), \partial_t v_\nu (t,\cdot)\rangle_{L^2}}\\[0.3 cm]
&\displaystyle{\ \ \ \ \ -(t+\tau)\langle c(t,\cdot)v_\nu(t,\cdot), \partial_t v_\nu (t,\cdot)\rangle_{L^2}}\\[0.3 cm]
&\displaystyle{\ \ \ \ \ \ +(t+\tau)\langle X_\nu(t,\cdot), \partial_t v_\nu (t,\cdot)\rangle_{L^2}.}
\end{array}
\end{equation}
Let now $\rho\in C^\infty_0({\mathbb R})$, with ${\rm supp}\, \rho\subseteq [-{1\over 2}, 
{1\over 2}]$, $\int_{\mathbb R}\rho(s)\, ds=1$ and $\rho(s)\geq 0$ for all $s\in{\mathbb R}$. We set, for $\varepsilon\in\,]0,1]$,
$$
a_{i,j,\varepsilon}(t,x):=\int_{\mathbb R} a_{i,j}(s,x){1\over  \varepsilon}\rho({t-s\over \varepsilon})\, ds.
$$
We deduce that, for all $\varepsilon\in\,]0,1]$,
\begin{equation}
\label{*}
k|\xi|^2\leq \sum_{i,j}a_{i,j,\varepsilon }(t,x)\xi_i\xi_j\leq k^{-1}|\xi|^2,
\end{equation}
\begin{equation}
\label{**}
|a_{i,j,\varepsilon}(t,x)-a_{i,j}(t,x)|\leq A_{LL}\mu(\varepsilon),
\end{equation}
and
\begin{equation}
\label{***}
|\partial_t a_{i,j,\varepsilon}(t,x)|\leq A_{LL}\|\rho'\|_{L^1}{\mu(\varepsilon)\over \varepsilon}.
\end{equation}
We set 
$$
a_{i,j,\nu}:= a_{i,j, \varepsilon}\quad {\rm with}\quad \varepsilon=2^{-2\nu},
$$ 
and 
$$
A_\nu:=(a_{i,j, \nu}(t,x))^n_{i,j=1}.
$$

In the second and fourth  term of the right hand side part of 
 (\ref{****}) we replace $A$ with $(A-A_\nu)+A_\nu$ and  $\partial_t v_\nu$ with the quantity given by  (\ref{eq2}), respectively.    We obtain

\begin{equation}
\begin{array}{ll}
\displaystyle{(t+\tau)\|\partial_t v_\nu(t, \cdot)\|^2_{L^2}}\\[0.3 cm]
\quad\displaystyle{={d\over dt}({\gamma\over 2}(t+\tau)\|v_\nu(t,\cdot)\|^2_{L^2})-{\gamma\over 2}\|v_\nu(t,\cdot)\|^2_{L^2}}\\[0.3 cm]
\displaystyle{\qquad\,+{d\over dt}({1\over 2}(t+\tau)\langle A_\nu(t,\cdot)\nabla_x v_\nu(t,\cdot), \nabla_x v_\nu(t,\cdot)\rangle_{(L^2)^n})}
\\[0.3 cm]
\displaystyle{\qquad\,\,-{1\over 2}\langle A_\nu(t,\cdot)\nabla_x v_\nu(t,\cdot), \nabla_x v_\nu(t,\cdot)\rangle_{(L^2)^n}}\\[0.3 cm]
\displaystyle{\qquad\,\,\,-{1\over 2}(t+\tau)\langle \partial_t A_\nu(t,\cdot)\nabla_x v_\nu(t,\cdot), \nabla_x v_\nu(t,\cdot)\rangle_{(L^2)^n}}\\[0.3 cm]
\displaystyle{\qquad\,\,\,\,+(t+\tau)\langle (A(t,\cdot)-A_\nu(t,\cdot))\nabla_x v_\nu(t,\cdot), \nabla_x\partial_t v_\nu(t,\cdot)\rangle_{(L^2)^n}}\\[0.3 cm]
\displaystyle{\qquad\,\,\,\,\,-{d\over dt}({1\over 2}(t+\tau)\Phi_\lambda'({t+\tau\over\beta})\|v_\nu(t,\cdot)\|^2_{L^2})
+{1\over 2}\Phi_\lambda'({t+\tau\over\beta})\|v_\nu(t,\cdot)\|^2_{L^2}}\\[0.3 cm]
\displaystyle{\qquad\,\,\,\,\,\,+{1\over 2}{t+\tau\over\beta} \Phi_\lambda''({t+\tau\over\beta})\|v_\nu(t,\cdot)\|^2_{L^2}
-\alpha\gamma (\log 2)\nu(t+\tau)\|v_\nu(t,\cdot)\|^2_{L^2}} \\[0.3 cm]
\displaystyle{\qquad\,\,\,\,\,\,\,-\alpha (\log 2)\nu(t+\tau)\langle A(t,\cdot)\nabla_x v_\nu(t,\cdot), \nabla_x v_\nu(t,\cdot)\rangle_{(L^2)^n}}\\[0.3 cm]
\displaystyle{\qquad\,\,\,\,\,\,\,\,+\alpha (\log 2)\nu(t+\tau)\Phi_\lambda'({t+\tau\over\beta})\|v_\nu(t,\cdot)\|^2_{L^2}}
\\[0.3 cm]\displaystyle{\qquad\,\,\,\,\,\,\,\,\,
+\alpha^2 (\log 2)^2\nu^2(t+\tau)\|v_\nu(t,\cdot)\|^2_{L^2}}\\[0.3 cm]
\qquad\displaystyle{\,\,\,\,\,\,\,\,\,\,+\alpha (\log 2)\nu(t+\tau)\langle v_\nu(t,\cdot), B(t,\cdot)\cdot\nabla_x v_\nu(t,\cdot)
\rangle_{L^2}}\\[0.3 cm]
\qquad\displaystyle{\,\,\,\,\,\,\,\,\,\,\,+\alpha (\log 2)\nu(t+\tau)\langle v_\nu(t,\cdot), c(t,\cdot) v_\nu(t,\cdot)
\rangle_{L^2}}
\\[0.3 cm]\,\,\,\,\,\,\,\,\,\,\,\,
\displaystyle{\qquad-\alpha (\log 2)\nu(t+\tau)\langle v_\nu(t,\cdot), X_\nu(t,\cdot) 
\rangle_{L^2}}
\\[0.3 cm]\displaystyle{\qquad\,\,\,\,\,\,\,\,\,\,\,\,\,
-(t+\tau)\langle B(t,\cdot)\cdot\nabla_x v_\nu(t,\cdot), \partial_t v_\nu(t,\cdot) 
\rangle_{L^2}}
\\[0.3 cm]
\displaystyle{\qquad\,\,\,\,\,\,\,\,\,\,\,\,\,\,\,-(t+\tau)\langle c(t,\cdot) v_\nu(t,\cdot), \partial_t v_\nu(t,\cdot) 
\rangle_{L^2}}
\\[0.3 cm]\displaystyle{\qquad\,\,\,\,\,\,\,\,\,\,\,\,\,\,\,\,
+(t+\tau)\langle  X_\nu(t,\cdot), 
 \partial_t v_\nu(t,\cdot)\rangle_{L^2}.}
\end{array}
\label{eq3}
\end{equation}

\subsection{Estimate for $\nu=0$}

We consider (\ref{eq3}) in the case of $\nu=0$. Using H\"older inequality, the inequalities (\ref{**}), (\ref{***}) and the fact that $\|\nabla_x v_0\|_{(L^2)^n}\leq 2\|v_0\|_{L^2}$ and 
the similar inequality for $\partial_t v_0$,  we deduce that, for $t\in [0,\sigma]$, 
$$
\begin{array}{ll}
\displaystyle{(t+\tau)\|\partial_t v_0(t, \cdot)\|^2_{L^2}}\\[0.3 cm]
\displaystyle{\,\quad\leq {d\over dt}({\gamma\over 2}(t+\tau)\|v_0(t,\cdot)\|^2_{L^2})-{\gamma\over 2}\|v_0(t,\cdot)\|^2_{L^2}}\\[0.3 cm]
\displaystyle{\,\qquad+{d\over dt}({1\over 2}(t+\tau)\langle A_0(t,\cdot)\nabla_x v_0(t,\cdot), \nabla_x v_0(t,\cdot)\rangle_{(L^2)^n})}
\\[0.3 cm]
\displaystyle{\,\,\qquad-{1\over 2}\langle A_0(t,\cdot)\nabla_x v_0(t,\cdot), \nabla_x v_0(t,\cdot)\rangle_{(L^2)^n}}
\\[0.3 cm]
\displaystyle{\qquad\,\,\,
+2n A_{LL}\|\rho'\|_{L^1}(t+\tau)\|v_0(t,\cdot)\|^2_{L^2}}\\[0.3 cm]
\displaystyle{\,\,\,\,\qquad+32n^2 A^2_{LL}(t+\tau)\|v_0(t,\cdot)\|^2_{L^2}+{1\over 8}(t+\tau)\|\partial_t v_0(t,\cdot)\|^2_{L^2}}\\[0.3 cm]
\displaystyle{\,\,\,\,\,\qquad-{d\over dt}({1\over 2}(t+\tau)\Phi_\lambda'({t+\tau\over\beta})\|v_0(t,\cdot)\|^2_{L^2})}
\\[0.3 cm]
\displaystyle{\qquad\,\,\,\,\,\,
+{1\over 2}\Phi_\lambda'({t+\tau\over\beta})\|v_0(t,\cdot)\|^2_{L^2}
+{1\over 2}{t+\tau\over\beta} \Phi_\lambda''({t+\tau\over\beta})\|v_0(t,\cdot)\|^2_{L^2}}
\\[0.3 cm]\displaystyle{\qquad\,\,\,\,\,\,\,
+8n B^2(t+\tau)\|v_0(t,\cdot)\|^2_{L^2}+{1\over 8}(t+\tau)\|\partial_t v_0(t,\cdot)\|^2_{L^2}}\\[0.3 cm]\displaystyle{\qquad\,\,\,\,\,\,\,\,
+2 C^2(t+\tau)\|v_0(t,\cdot)\|^2_{L^2}+{1\over 8}(t+\tau)\|\partial_t v_0(t,\cdot)\|^2_{L^2}}
\\[0.3 cm]\displaystyle{\qquad\,\,\,\,\,\,\,\,\,
+(t+\tau)\langle  X_\nu(t,\cdot) ,
 \partial_t v_0(t,\cdot)\rangle_{L^2}.}
\end{array}
$$

Choosing $\gamma$ such that ${\gamma \over 4}\geq (2n A_{LL}\|\rho'\|_{L^1} +32n^2 A^2_{LL}+
8n B^2+ 2C^2)(\sigma+\tau)$
the term
$$
\begin{array}{ll} 
+2n A_{LL}\|\rho'\|_{L^1}(t+\tau)\|v_0(t,\cdot)\|^2_{L^2}+32n^2 A^2_{LL}(t+\tau)\|v_0(t,\cdot)\|^2_{L^2}\\[0.2 cm]
\qquad+8n B^2(t+\tau)\|v_0(t,\cdot)\|^2_{L^2}+2 C^2(t+\tau)\|v_0(t,\cdot)\|^2_{L^2}
\end{array}
$$
is absorbed by $-{\gamma\over 4}\|v_0(t,\cdot)\|^2_{L^2}$.

Recalling that $\Phi_\lambda$ satisfies (\ref{ode}), i. e.
$$
y\, \Phi_\lambda''(y) =-\lambda\, \Phi_\lambda'(y) (1+|\log {1\over \Phi_\lambda'(y)}|),
$$
with $\lambda>1$, the term ${1\over 2}\Phi_\lambda'({t+\tau\over\beta})\|v_0(t,\cdot)\|^2_{L^2}$ is balanced by ${1\over 2}{t+\tau\over\beta} \Phi_\lambda''({t+\tau\over\beta})\|v_0(t,\cdot)\|^2_{L^2}$.

We obtain
$$
\begin{array}{ll}
\displaystyle{
{5\over 8}(t+\tau)\|\partial_t v_0(t, \cdot)\|^2_{L^2}\leq}&\displaystyle{ {d\over dt}({\gamma\over 2}(t+\tau)\|v_0(t,\cdot)\|^2_{L^2})-{\gamma\over 4}\|v_0(t,\cdot)\|^2_{L^2}
}\\[0.3 cm]
&\displaystyle{+{d\over dt}({1\over 2}(t+\tau)\langle A_0(t,\cdot)\nabla_x v_0(t,\cdot), \nabla_x v_0(t,\cdot)\rangle_{(L^2)^n})}
\\[0.3 cm]&\displaystyle{\quad
-{1\over 2}\langle A_0(t,\cdot)\nabla_x v_0(t,\cdot), \nabla_x v_0(t,\cdot)\rangle_{(L^2)^n}}
\\[0.3 cm]&\displaystyle{\qquad
-{d\over dt}({1\over 2}(t+\tau)\Phi_\lambda'({t+\tau\over\beta})\|v_0(t,\cdot)\|^2_{L^2})}
\\[0.3 cm]&\displaystyle{\qquad\quad
+(t+\tau)\langle  X_\nu(t,\cdot) ,
 \partial_t v_0(t,\cdot)\rangle_{L^2}.}
\end{array}
$$
Integrating the previous inequality between $0$ and $s$, with $s\leq \sigma $, we have
$$
\begin{array}{ll}
\displaystyle{{1\over 2}\int_0^s \langle A_0(t,\cdot)\nabla_x v_0(t,\cdot), \nabla_x
v_0(t,\cdot)\rangle_{(L^2)^n}\, dt +{\gamma\over 8}\int_0^s \|v_0(t,\cdot)\|^2_{L^2}\, dt}
\\[0.3 cm]
\qquad\qquad\displaystyle{\leq{1\over 2}(s+\tau)\langle A_0(s,\cdot)\nabla_x v_0(s,\cdot), \nabla_x v_0(s,\cdot)\rangle_{(L^2)^n}}
\\[0.3 cm]
\quad\qquad\qquad\displaystyle{+{\gamma\over 2}(s+\tau)\|v_0(s,\cdot)\|^2_{L^2}
+{1\over 2}\tau\,\Phi_\lambda'({\tau\over\beta})\|v_0(0,\cdot)\|^2_{L^2}}
\\[0.3 cm]
\qquad\qquad\qquad\displaystyle{-{\gamma\over 8}\int_0^s \|v_0(t,\cdot)\|^2_{L^2}\, dt
-{5\over 8}\int_0^s (t+\tau)\|\partial_t v_0(t,\cdot)\|^2_{L^2}\, dt}
\\[0.3 cm]
\quad\qquad\qquad\qquad\displaystyle{+\int_0^s (t+\tau)\langle  X_\nu(t,\cdot), 
 \partial_t v_0(t,\cdot)\rangle_{L^2}\, dt,}
\end{array}
$$
where on the right hand side part some negative terms have been neglected. Again from the fact that $\|\nabla_x v_0\|_{(L^2)^n}\leq 2\|v_0\|_{L^2}$, using also (\ref{*}), we finally deduce
\begin{equation}
\begin{array}{ll}
\displaystyle{{\gamma\over 8}\int_0^s \|v_0(t,\cdot)\|^2_{L^2}\, dt}&\!\!\!\leq
\displaystyle{(2k^{-1} +{\gamma \over 2})(s+\tau)\|v_0(s,\cdot)\|^2_{L^2}}
\\[0.3 cm]&\ 
\displaystyle{+{1\over 2}\tau\,\Phi_\lambda'({\tau\over\beta})\|v_0(0,\cdot)\|^2_{L^2}
-{\gamma\over 8}\int_0^s \|v_0(t,\cdot)\|^2_{L^2}\, dt}
\\[0.3 cm]&\quad\ 
\displaystyle{
-{5\over 8}\int_0^s  (t+\tau) \|\partial_t v_0(t,\cdot)\|^2_{L^2}\, dt}
\\[0.3 cm]&\qquad\ 
\displaystyle{+\int_0^s (t+\tau)\langle  X_\nu(t,\cdot), 
 \partial_t v_0(t,\cdot)\rangle_{L^2}\, dt.}
\end{array}
\label{stella}
\end{equation}

\subsection{Estimate for  $\nu\geq 1$}

We consider (\ref{eq3}) in the case of $\nu\geq 1$. Using again  H\"older inequality, the inequalities (\ref{**}), (\ref{***}) and (\ref{d1}),  we have that
$$
\begin{array}{ll}\displaystyle{
(t+\tau)\|\partial_t v_\nu(t, \cdot)\|^2_{L^2}}\\[0.3 cm]
\quad\displaystyle{\leq{d\over dt}({\gamma\over 2}(t+\tau)\|v_\nu(t,\cdot)\|^2_{L^2})-{\gamma\over 2}\|v_\nu(t,\cdot)\|^2_{L^2}}\\[0.3 cm]
\displaystyle{\qquad\,+{d\over dt}({1\over 2}(t+\tau)\langle A_\nu(t,\cdot)\nabla_x v_\nu(t,\cdot), \nabla_x v_\nu(t,\cdot)\rangle_{(L^2)^n})}
\\[0.3 cm]\displaystyle{
\qquad\,\,-{1\over 2}\langle A_\nu(t,\cdot)\nabla_x v_\nu(t,\cdot), \nabla_x v_\nu(t,\cdot)\rangle_{(L^2)^n}}\\[0.3 cm]
\displaystyle{\qquad\,\,\,+2(1+2\log 2)n A_{LL}\|\rho'\|_{L^1}\nu(t+\tau) 2^{2\nu}\|v_\nu(t,\cdot)\|^2_{L^2}}\\[0.3 cm]
\displaystyle{\qquad\,\,\,\,+32(1+2\log 2)^2n^2 A_{LL}^2\nu^2(t+\tau) \|v_\nu(t,\cdot)\|^2_{L^2}+{1\over 8} (t+\tau)\|\partial_t v_\nu(t,\cdot)\|^2_{L^2}}
\\[0.3 cm]
\displaystyle{\qquad\,\,\,\,\,-{d\over dt}({1\over 2}(t+\tau)\Phi_\lambda'({t+\tau\over\beta})\|v_\nu(t,\cdot)\|^2_{L^2})
+{1\over 2}\Phi_\lambda'({t+\tau\over\beta})\|v_\nu(t,\cdot)\|^2_{L^2}}\\[0.3 cm]
\displaystyle{\qquad\,\,\,\,\,\,+{1\over 2}{t+\tau\over\beta} \Phi_\lambda''({t+\tau\over\beta})\|v_\nu(t,\cdot)\|^2_{L^2}
-\alpha\gamma (\log 2)\nu(t+\tau)\|v_\nu(t,\cdot)\|^2_{L^2}} \\[0.3 cm]
\displaystyle{\qquad\,\,\,\,\,\,\,-\alpha (\log 2){k\over 4}\nu(t+\tau)2^{2\nu}\|v_\nu(t,\cdot)\|^2_{L^2}
+\alpha (\log 2)\nu(t+\tau)\Phi_\lambda'({t+\tau\over\beta})\|v_\nu(t,\cdot)\|^2_{L^2}}
\\[0.3 cm]\displaystyle{\qquad\,\,\,\,\,\,\,\,
+\alpha^2 (\log 2)^2\nu^2(t+\tau)\|v_\nu(t,\cdot)\|^2_{L^2}}
\\[0.3 cm]\displaystyle{\qquad\,\,\,\,\,\,\,\,\,
+\alpha 2(\log 2)n^{1/2}B \nu(t+\tau)2^\nu \|v_\nu(t,\cdot)\|^2_{L^2}}
\\[0.3 cm]\displaystyle{\qquad\,\,\,\,\,\,\,\,\,\,
+\alpha(\log 2) C\nu(t+\tau)\|v_\nu(t,\cdot)\|^2_{L^2}}
\\[0.3 cm]\displaystyle{\qquad\,\,\,\,\,\,\,\,\,\,\,
+32n B^2 (t+\tau)2^{2\nu}\|v_\nu(t,\cdot)\|^2_{L^2}
+{1\over 8} (t+\tau)\|\partial_t v_\nu(t,\cdot)\|^2_{L^2}}
\\[0.3 cm]
\displaystyle{\qquad\,\,\,\,\,\,\,\,\,\,\,\,+2C^2(t+\tau)\|v_\nu(t,\cdot)\|^2_{L^2}
+{1\over 8} (t+\tau)\|\partial_t v_\nu(t,\cdot)\|^2_{L^2}}
\\[0.3 cm]
\displaystyle{\,\,\,\,\,\,\,\,\,\,\,\,\,
\qquad-\alpha (\log 2)\nu(t+\tau)\langle v_\nu(t,\cdot), X_\nu(t,\cdot) 
\rangle_{L^2}}
\\[0.3 cm]\displaystyle{\qquad\,\,\,\,\,\,\,\,\,\,\,\,\,\,\,
+(t+\tau)\langle  X_\nu(t,\cdot) ,
 \partial_t v_\nu(t,\cdot)\rangle_{L^2}.}
\end{array}
$$

Let now $\alpha=\max\{T^{-1}, \alpha_1\}$, where
$$
\alpha_1:={16\over k\log 2}(2(1+2\log 2)n A_{LL}\|\rho'\|_{L^1}
+32(1+2\log 2)^2n^2 A_{LL}^2+32n B^2),
$$
then
$$
\begin{array}{ll}
\displaystyle{-{\alpha\over 4}(\log 2){k\over 4}\nu 2^{2\nu}+
2(1+2\log 2)n A_{LL}\|\rho'\|_{L^1}\nu 2^{2\nu}}\\[0.3 cm]
\qquad\qquad\qquad\qquad\qquad \displaystyle{+32(1+2\log 2)^2n^2 A_{LL}^2\nu^2+32n B^22^{2\nu}\leq 0,}
\end{array}
$$
and the term 
$$
\begin{array}{ll}
\displaystyle{2(1+2\log 2)n A_{LL}\|\rho'\|_{L^1}\nu(t+\tau) 2^{2\nu}\|v_\nu(t,\cdot)\|^2_{L^2}}\\[0.3 cm]
\displaystyle{\qquad+32(1+2\log 2)^2n^2 A_{LL}^2\nu^2(t+\tau) \|v_\nu(t,\cdot)\|^2_{L^2}}
\\[0.2 cm]\displaystyle{\qquad\qquad
+32n B^2 (t+\tau)2^{2\nu}\|v_\nu(t,\cdot)\|^2_{L^2}}
\end{array}
$$
is absorbed by $-{\alpha\over 4} (\log 2){k\over 4}\nu(t+\tau)2^{2\nu}\|v_\nu(t,\cdot)\|^2_{L^2}. $

Since 
$y\, \Phi_\lambda''(y) =-\lambda\, \Phi_\lambda'(y) (1+|\log {1\over \Phi_\lambda'(y)}|)$, supposing  $\lambda>2$ we have
$$
{1\over 4}{t+\tau\over \beta}\Phi_\lambda''({t+\tau\over \beta})
\leq -{1\over 2}\Phi_\lambda'({t+\tau\over \beta}),
$$
and the term ${1\over 2}\Phi_\lambda'({t+\tau\over\beta})\|v_\nu(t,\cdot)\|^2_{L^2}$ is absorbed by ${1\over 4}{t+\tau\over\beta} \Phi_\lambda''({t+\tau\over\beta})\|v_\nu(t,\cdot)\|^2_{L^2}$.

Consider now the term 
$$
\alpha (\log 2)\nu(t+\tau)\Phi_\lambda'({t+\tau\over\beta})\|v_\nu(t,\cdot)\|^2_{L^2}.
$$
Let $k'=\min\{k, 16\}$. If $\nu\geq (\log 2)^{-1}\log ({16\over k'}\Phi_\lambda'({t+\tau\over\beta}))$, then
$$
-{\alpha\over 4}(\log 2){k\over 4}\nu 2^{2\nu}\leq -\alpha(\log 2)\Phi_\lambda'({t+\tau\over\beta})\nu.
$$
On the contrary, if $\nu< (\log 2)^{-1}\log ({16\over k'}\Phi_\lambda'({t+\tau\over\beta}))$
then ${16\over k'}\Phi_\lambda'({t+\tau\over\beta})>2^\nu$, so that
$$
\begin{array}{ll}
\displaystyle{{1\over 4}{t+\tau\over\beta}\Phi_\lambda''({t+\tau\over\beta})}&=
\displaystyle{-{1\over 4}\lambda( \Phi_\lambda'({t+\tau\over\beta}))^2
\mu({1\over \Phi_\lambda'({t+\tau\over\beta})})}\\[0.4 cm]
& \leq\displaystyle{-{1\over 4}\lambda( \Phi_\lambda'({t+\tau\over\beta}))^2
\mu({1\over {16\over k'}\Phi_\lambda'({t+\tau\over\beta})})}\\[0.4 cm]
& \leq\displaystyle{-{1\over 4}\lambda {k'\over 16} \Phi_\lambda'({t+\tau\over\beta})
(1+|\log({1\over {16\over k'}\Phi_\lambda'({t+\tau\over\beta})})|)}\\[0.3 cm]
& \leq\displaystyle{-{1\over 4}\lambda {k'\over 16} \Phi_\lambda'({t+\tau\over\beta})
(1+\nu \log 2)}\\[0.3 cm]
& \leq \displaystyle{-\lambda {k'\log 2\over 48} \Phi_\lambda'({t+\tau\over\beta})
\nu},
\end{array}
$$
where we have used the fact that the function $\mu$ is increasing.
Consequently if we choose $\lambda$ in such a way that
$$
{\lambda k'\log 2\over 48}\geq \alpha (\log 2)(\sigma+\tau), 
\quad{\rm i.\ e.}\quad \lambda\geq{48\alpha(\sigma+\tau)\over k'} ,
$$
then, if $\nu< (\log 2)^{-1}\log ({16\over k'}\Phi_\lambda'({t+\tau\over\beta}))$, we have
$$
{1\over 4}{t+\tau\over\beta}\Phi_\lambda''({t+\tau\over\beta})\leq -\alpha(\log 2)({t+\tau})\Phi_\lambda'({t+\tau\over\beta})\nu.
$$
In conclusion $\alpha (\log 2)\nu(t+\tau)\Phi_\lambda'({t+\tau\over\beta})\|v_\nu(t,\cdot)\|^2_{L^2}$ is balanced by $-{\alpha\over 4} (\log 2){k\over 4}\nu(t+\tau)2^{2\nu}\|v_\nu(t,\cdot)\|^2_{L^2}+{1\over 4}{t+\tau\over\beta} \Phi_\lambda''({t+\tau\over\beta})\|v_\nu(t,\cdot)\|^2_{L^2}$. We remark that the computations here above are the main
 cause for the introduction of the weight function $ \Phi_\lambda$.

Consider now the sum
$$
\alpha^2 (\log 2)^2(t+\tau)\nu^2\|v_\nu(t,\cdot)\|^2_{L^2}
+\alpha 2(\log 2)n^{1/2}B (t+\tau)\nu2^\nu \|v_\nu(t,\cdot)\|^2_{L^2}.
$$
If $\nu\geq (\log 2)^{-1}\log ({1\over k}(16\alpha \log 2+32 n^{1/2}B))=:\bar \nu_1$ then 
$$
-{\alpha\over 4}(\log 2){k\over 4}\nu 2^{2\nu}+(\alpha^2 (\log 2)^2\nu^2
+\alpha 2(\log 2)n^{1/2}B \nu2^\nu)\leq 0.
$$
If $\nu< (\log 2)^{-1}\log ({1\over k}(16\alpha \log 2+32 n^{1/2}B))=\bar \nu_1$, choosing $\gamma$ in such a way that
$$
{\gamma\over 4}\geq (\alpha^2 (\log 2)^2\bar\nu_1^2
+\alpha 2(\log 2)n^{1/2}B \bar\nu_12^{\bar\nu_1})(\sigma+\tau),
$$ 
we obtain 
$$
-{\gamma\over 4}+(\alpha^2 (\log 2)^2\nu^2
+\alpha 2(\log 2)n^{1/2}B \nu 2^{\nu})(t+\tau)\leq 0,
$$
and consequently the term $\alpha^2 (\log 2)^2(t+\tau)\nu^2\|v_\nu(t,\cdot)\|^2_{L^2}
+\alpha 2(\log 2)n^{1/2}B (t+\tau)\nu2^\nu \|v_\nu(t,\cdot)\|^2_{L^2}$ 
is absorbed by $-{\alpha\over 4} (\log 2){k\over 4}\nu(t+\tau)2^{2\nu}\|v_\nu(t,\cdot)\|^2_{L^2}-{\gamma\over 4}\|v_\nu(t,\cdot)\|^2_{L^2}$.

Consider finally
$$
\alpha(\log 2) C(t+\tau)\nu\|v_\nu(t,\cdot)\|^2_{L^2}
+2C^2(t+\tau)\|v_\nu(t,\cdot)\|^2_{L^2}.
$$
If we take $\gamma$ such that
$$
\gamma\geq {\alpha(\log 2)C+2C^2\over \alpha\log 2},
$$
then 
$$
-\alpha\gamma(\log 2)\nu+\alpha(\log 2)C\nu+2C^2\leq 0,
$$
and the above quoted term is absorbed by $-\alpha \gamma (\log 2) \nu  \|v_\nu(t,\cdot)\|^2_{L^2}$.

Summing up, we set 
$$
\alpha_1:={16\over k\log 2}(2(1+2\log 2)n A_{LL}\|\rho'\|_{L^1}
+32(1+2\log 2)^2n^2 A_{LL}^2+32n B^2),
$$
$$
\alpha:=\max\{T^{-1}, \alpha_1\}, \qquad \sigma:=1/\alpha,\qquad k':=\min\{k, 16\}.
$$
We choose $\tau\in \,]0,\sigma/2[$ and we define
$$
\bar \nu_1:=(\log 2)^{-1}\log ({1\over k}(16\alpha \log 2+32 n^{1/2}B)).
$$
We choose $\lambda$, $\gamma$ in such a way that
$$
\lambda\geq \max\{2, {48\alpha(\sigma+\tau)\over k'}\},
$$
$$
\gamma\geq \max\{{\alpha(\log 2)C+2C^2\over \alpha\log 2},
4(\alpha^2 (\log 2)^2\bar\nu_1^2
+\alpha 2(\log 2)n^{1/2}B \bar\nu_12^{\bar\nu_1})(\sigma+\tau)\}.
$$
Then, for all $\beta\geq \sigma+\tau$ and for all $\nu\geq 1$ we have
$$
\begin{array}{ll}\displaystyle{
{5\over 8}(t+\tau)\|\partial_t v_\nu(t, \cdot)\|^2_{L^2}\leq}&\displaystyle{ {d\over dt}({\gamma\over 2}(t+\tau)\|v_\nu(t,\cdot)\|^2_{L^2})-{\gamma\over 4}\|v_\nu(t,\cdot)\|^2_{L^2}}
\\[0.3 cm]
&\displaystyle{+{d\over dt}({1\over 2}(t+\tau)\langle A_\nu(t,\cdot)\nabla_x v_\nu(t,\cdot), \nabla_x v_\nu(t,\cdot)\rangle_{(L^2)^n})}
\\[0.3 cm]&\displaystyle{\quad
-{1\over 2}\langle A_\nu(t,\cdot)\nabla_x v_\nu(t,\cdot), \nabla_x v_\nu(t,\cdot)\rangle_{(L^2)^n}}
\\[0.3 cm]&\displaystyle{\qquad
-{d\over dt}({1\over 2}(t+\tau)\Phi_\lambda'({t+\tau\over\beta})\|v_\nu(t,\cdot)\|^2_{L^2})}
\\[0.3 cm]&\displaystyle{\qquad\quad
-{\alpha\over 4}(\log 2){k\over 4}(t+\tau)\nu 2^\nu\|v_\nu(t,\cdot)\|^2_{L^2}}
\\[0.3 cm]&\displaystyle{\qquad\qquad
-\alpha(\log 2)(t+\tau)\nu \langle v_\nu(t,\cdot), X_\nu(t,\cdot)\rangle}
\\[0.3 cm]&\displaystyle{\qquad\qquad\quad
+(t+\tau)\langle  X_\nu(t,\cdot) ,
 \partial_t v_0(t,\cdot)\rangle_{L^2}.}
\end{array}
$$
We integrate this last  inequality between $0$ and $s$, with $s\leq \sigma $, and we obtain
$$
\begin{array}{ll}
\displaystyle{{1\over 2}\int_0^s \langle A_\nu(t,\cdot)\nabla_x v_\nu(t,\cdot), \nabla_x
v_\nu(t,\cdot)\rangle_{(L^2)^n}\, dt +{\gamma\over 8}\int_0^s \|v_\nu(t,\cdot)\|^2_{L^2}\, dt}
\\[0.3 cm]
\qquad\qquad\qquad\displaystyle{\leq{1\over 2}(s+\tau)\langle A_\nu(s,\cdot)\nabla_x v_\nu(s,\cdot), \nabla_x v_\nu(s,\cdot)\rangle_{(L^2)^n}}
\\[0.3 cm]
\quad \qquad\qquad\qquad\displaystyle{+{\gamma\over 2}(s+\tau)\|v_\nu(s,\cdot)\|^2_{L^2}
+{1\over 2}\tau\,\Phi_\lambda'({\tau\over\beta})\|v_\nu(0,\cdot)\|^2_{L^2}}
\\[0.3 cm]
\quad \quad \qquad\qquad\qquad\displaystyle{-{\gamma\over 8}\int_0^s \|v_\nu(t,\cdot)\|^2_{L^2}\, dt
-{5\over 8}\int_0^s (t+\tau)\|\partial_t v_\nu(t,\cdot)\|^2_{L^2}\, dt}
\\[0.3 cm]
\quad \qquad\qquad\qquad\qquad\displaystyle{-{\alpha k (\log 2)\over 16} \int_0^s(t+\tau)\nu 2^{2\nu} \|v_\nu(t,\cdot)\|^2_{L^2}\, dt}
\\[0.3 cm]
\qquad\qquad \qquad\qquad\qquad\displaystyle{-\alpha(\log 2) \int_0^s(t+\tau)\nu \langle v_\nu(t,\cdot), X_\nu \rangle_{L^2}\, dt,}
\\[0.3 cm]
\qquad\qquad\quad\qquad \qquad\qquad\displaystyle{+\int_0^s (t+\tau)\langle  X_\nu(t,\cdot), 
 \partial_t v_\nu(t,\cdot)\rangle_{L^2}\, dt,}
\end{array}
$$
where again on the right hand side part some negative terms have been neglected. Using  (\ref{*}) and (\ref{d1}) we  obtain
\begin{equation}
\begin{array}{ll}
\displaystyle{{k\over 8}\int_0^s 2^{2\nu}\|v_\nu(t,\cdot)\|^2_{L^2}\, dt+{\gamma\over 8}\int_0^s \|v_\nu(t,\cdot)\|^2_{L^2}\, dt}\\[0.3 cm]\qquad\qquad
\displaystyle{ \leq(2^{2\nu+1}k^{-1} +{\gamma \over 2})(s+\tau)\|v_\nu(s,\cdot)\|^2_{L^2}
+{1\over 2}\tau\,\Phi_\lambda'({\tau\over\beta})\|v_\nu(0,\cdot)\|^2_{L^2}}
\\[0.3 cm]\qquad\qquad\quad
\displaystyle{-{\gamma\over 8}\int_0^s \|v_\nu(t,\cdot)\|^2_{L^2}\, dt
-{5\over 8}\int_0^s (t+\tau)\|\partial_t v_\nu(t,\cdot)\|^2_{L^2}\, dt}
\\[0.3 cm]\qquad\qquad\qquad
\displaystyle{-{\alpha k (\log 2)\over 16} \int_0^s(t+\tau)\nu 2^{2\nu} \|v_\nu(t,\cdot)\|^2_{L^2}\, dt}
\\[0.3 cm]\qquad\qquad\qquad\quad
\displaystyle{-\alpha(\log 2) \int_0^s(t+\tau)\nu \langle v_\nu(t,\cdot), X_\nu \rangle_{L^2}\, dt,}
\\[0.3 cm]\qquad\qquad\qquad\qquad
\displaystyle{+\int_0^s (t+\tau)\langle  X_\nu(t,\cdot), 
 \partial_t v_0(t,\cdot)\rangle_{L^2}\, dt.}
\end{array}
\label{stella2}
\end{equation}

\subsection{Estimate for the commutator term}

We collect together (\ref{stella}) and (\ref{stella2}). We deduce that
$$
\begin{array}{ll}
\displaystyle{{k\over 8}\int_0^s \sum_{\nu=1}^{+\infty} 2^{2\nu}\|v_\nu(t,\cdot)\|^2_{L^2}\, dt+{\gamma\over 8}\int_0^s \sum_{\nu=0}^{+\infty} \|v_\nu(t,\cdot)\|^2_{L^2}\, dt}\\[0.3 cm]\qquad\qquad
\displaystyle{\leq  {2\over k}(s+\tau)\sum_{\nu=0}^{+\infty}2^{2\nu}\|v_\nu(s,\cdot)\|^2_{L^2}
+{\gamma \over 2}(s+\tau)\sum_{\nu=0}^{+\infty}\|v_\nu(s,\cdot)\|^2_{L^2}}
\\[0.3 cm]\qquad\qquad\quad
\displaystyle{+{1\over 2}\tau\,\Phi_\lambda'({\tau\over\beta})\sum_{\nu=0}^{+\infty}\|v_\nu(0,\cdot)\|^2_{L^2}}
\\[0.3 cm]\qquad\qquad\qquad
\displaystyle{-{\gamma\over 8}\int_0^s \sum_{\nu=0}^{+\infty}\|v_\nu(t,\cdot)\|^2_{L^2}\, dt -{5\over 8}\int_0^s (t+\tau)\sum_{\nu=0}^{+\infty}\|\partial_t v_\nu(t,\cdot)\|^2_{L^2}\, dt}
\\[0.3 cm]\qquad\qquad\qquad\quad
\displaystyle{-{\alpha k (\log 2)\over 16} \int_0^s(t+\tau)\sum_{\nu=0}^{+\infty}\nu 2^{2\nu} \|v_\nu(t,\cdot)\|^2_{L^2}\, dt}
\\[0.3 cm]\qquad\qquad\qquad\qquad
\displaystyle{-\alpha(\log 2) \int_0^s(t+\tau)\sum_{\nu=0}^{+\infty}\nu \langle v_\nu(t,\cdot), X_\nu \rangle_{L^2}\, dt,}
\\[0.3 cm]\qquad\qquad\qquad\qquad\quad
\displaystyle{+\int_0^s (t+\tau)\sum_{\nu=0}^{+\infty}\langle  X_\nu(t,\cdot), 
 \partial_t v_0(t,\cdot)\rangle_{L^2}\, dt.}
\end{array}
$$
Now we want  to estimate the last two terms: to do this we shall follow essentially the ideas contained in \cite{CL}. We recall that
$$
X_\nu(t,x)=X_\nu^1(t,x)+X_\nu^2(t,x)+X_\nu^3(t,x),
$$
where
$$
\begin{array}{ll}
\displaystyle{X^1_\nu(t,x):= -\sum_{i,j} \partial_{x_i}([\varphi_\nu(D_x),a_{i,j}(t,x)]2^{-\alpha \nu t}\partial_{x_j} w)},\\[0.3 cm]
\displaystyle{X^2_\nu(t,x):= -\sum_{j} [\varphi_\nu(D_x),b_j(t,x)]2^{-\alpha \nu t}\partial_{x_j} w},\\[0.3 cm]
\displaystyle{X^3_\nu(t,x):=-[\varphi_\nu(D_x),c(t,x)]2^{-\alpha \nu t}w.}
\end{array}
$$
We start with 
$$
\begin{array}{ll}
\displaystyle{\sum_{\nu=0}^{+\infty}\langle X^1_\nu(t,\cdot), \partial_t v_\nu(t, \cdot)\rangle_{L^2}}\\[0.3 cm]
\qquad\qquad\displaystyle{=\sum_{i,j}\sum_\nu \langle [\varphi_\nu(D_x),a_{i,j}]2^{-\alpha \nu t}\partial_{x_j} w, \partial_{x_i}\partial_t v_\nu\rangle_{L^2}}\\[0.3 cm]
\qquad\qquad\displaystyle{=\sum_{i,j}\sum_{\nu, \mu}
\langle [\varphi_\nu(D_x),a_{i,j}]2^{-\alpha \nu t}\partial_{x_j} w_\mu, \partial_{x_i}\partial_t v_\nu\rangle_{L^2}}\\[0.3 cm]
\qquad\qquad\displaystyle{=\sum_{i,j}\sum_{\nu, \mu}
\langle ([\varphi_\nu(D_x),a_{i,j}]\psi_\mu(D_x))(2^{-\alpha (\nu-\mu) t}\partial_{x_j} v_\mu), \partial_{x_i}\partial_t v_\nu\rangle_{L^2}}
\end{array}
$$
where $\psi_\mu(D_x):=\varphi_{\mu-1}(D_x)+\varphi_\mu(D_x)+\varphi_{\mu+1}(D_x)$. Consequently
$$
\begin{array}{ll}
\displaystyle{|\sum_{\nu=0}^{+\infty}\langle X^1_\nu(t,\cdot), \partial_t v_\nu(t, \cdot)\rangle_{L^2}|}\\[0.3 cm]
\ \ \displaystyle{\leq\sum_{i,j}\sum_{\nu, \mu}
\|[\varphi_\nu(D_x),a_{i,j}]\psi_\mu(D_x)\|_{{\mathcal L}(L^2,L^2)}2^{-\alpha (\nu-\mu) t}
2^{\mu+1}2^{\nu+1}\|v_\mu\|_{L^2}\|\partial_t v_\nu\|_{L^2}}.
\end{array}
$$
We know that there exists $Q_A>0$ such that
$$
\|[\varphi_\nu(D_x), a_{ij}]\psi_\mu(D_x)\|_{{\mathcal L}(L^2, L^2)}
\leq \left\{
\begin{array}{ll}
\displaystyle{3Q_A\, 2^{-2\nu}}&\quad{\rm if}\quad |\mu-\nu|\leq 2,\\[0.5 cm]
\displaystyle{3Q_A\, 2^{-2\max\{\nu, \mu\}}}\ &\quad{\rm if}\quad |\mu-\nu|\geq 3.
 \end{array}
 \right.
 $$
Setting $k_{\nu ,\mu}(t):=2^{-\alpha (\nu-\mu) t}
2^\nu \|[\varphi_\nu(D_x),a_{i,j}]\psi_\mu(D_x)\|_{{\mathcal L}(L^2,L^2)}$, for $0\leq t\leq 1/\alpha$, we get, for a fixed $\nu\geq 0$, 
$$
\begin{array}{ll}
\displaystyle{\sum_{\mu} |k_{\nu, \mu}(t)|}&\leq \displaystyle{\sum_{|\mu-\nu|\geq 3}
2^{-\alpha (\nu-\mu) t} 2^\nu 3Q_A 2^{-2\max\{\nu, \mu\}}}\\[0.3 cm]
&\qquad\qquad\qquad\qquad\qquad\quad\displaystyle{+\sum_{\mu=\nu-2}^{\mu=\nu+2}2^{-\alpha (\nu-\mu) t} 2^\nu 3Q_A 2^{-2\nu}}
\\[0.3 cm]
&\displaystyle{\leq\sum_{\mu=0}^{\nu-3}2^{-\alpha (\nu-\mu) t} 2^\nu 3Q_A 2^{-2\nu}
+\sum_{\mu=\nu+3}^{+\infty}2^{-\alpha (\nu-\mu) t} 2^\nu 3Q_A 2^{-2\mu}}\\[0.5 cm]
&\qquad\qquad\qquad\qquad\quad\displaystyle{+3Q_A(2^{-2\alpha t}+2^{-\alpha t}+1+2^{\alpha t}+2^{2\alpha t})}
\\[0.3 cm]
&\displaystyle{\leq 30\, Q_A.}
 \end{array}
 $$
On the other hand, for a fixed $\nu\geq 0$,
$$
\begin{array}{ll}
\displaystyle{\sum_{\nu} |k_{\nu, \mu}(t)|}&\leq \displaystyle{\sum_{|\mu-\nu|\geq 3}
2^{-\alpha (\nu-\mu) t} 2^\nu 3Q_A 2^{-2\max\{\nu, \mu\}}}\\[0.3 cm]
&\qquad\qquad\qquad\qquad\qquad\displaystyle{+\sum_{\nu=\mu-2}^{\nu=\mu+2}2^{-\alpha (\nu-\mu) t} 2^\nu 3Q_A 2^{-2\nu}}
\\[0.3 cm]
&\displaystyle{\leq\sum_{\nu=0}^{\mu-3}2^{-\alpha (\nu-\mu) t} 2^\nu 3Q_A 2^{-2\mu}
+\sum_{\nu=\mu+3}^{+\infty}2^{-\alpha (\nu-\mu) t} 2^\nu 3Q_A 2^{-2\nu}}\\[0.3 cm]
&\qquad\qquad\qquad\quad\displaystyle{+3Q_A(2^{2\alpha t}+2^{\alpha t}+1+2^{-\alpha t}+2^{-2\alpha t})}
\\[0.3 cm]
&\displaystyle{\leq 35\, Q_A.}
 \end{array}
 $$
From Schur's criterion it follows that
$$
|\sum_{\nu=0}^{+\infty}\langle X^1_\nu(t,\cdot), \partial_t v_\nu(t, \cdot)\rangle_{L^2}|
\leq n^2 \, 140\, Q_A(\sum_{\nu=0}^{+\infty}2^{2\mu}\|v_\mu\|^2_{L^2})^{1/2}
(\sum_{\nu=0}^{+\infty}\|\partial_t v_\mu\|^2_{L^2})^{1/2},
$$
and then, for all $\eta>0$, 
$$
|\sum_{\nu=0}^{+\infty}\langle X^1_\nu(t,\cdot), \partial_t v_\nu(t, \cdot)\rangle_{L^2}|
\leq {(n^2 \, 140\, Q_A)^2\over 2\eta}\sum_{\nu=0}^{+\infty}2^{2\mu}\|v_\mu\|^2_{L^2}+
{\eta\over 2} \sum_{\nu=0}^{+\infty}\|\partial_t v_\mu\|^2_{L^2}.
$$
Arguing in a similar way we deduce that for all $\eta>0$ there exists $Q_\eta>0$ such that
$$
|\sum_{\nu=0}^{+\infty}\langle X_\nu(t,\cdot), \partial_t v_\nu(t, \cdot)\rangle_{L^2}|
\leq Q_\eta \sum_{\nu=0}^{+\infty}2^{2\nu}\|v_\nu(t,\cdot)\|^2_{L^2}+
\eta \sum_{\nu=0}^{+\infty}\|\partial_t v_\nu(t,\cdot)\|^2_{L^2},
$$
and there exists $\tilde Q_1>0$ such that
$$
|\sum_{\nu=0}^{+\infty}\nu \langle v_\nu(t, \cdot),  X_\nu(t,\cdot) \rangle_{L^2}|
\leq \tilde Q_1 \sum_{\nu=0}^{+\infty} 2^{2\nu} \|v_\nu(t,\cdot)\|^2_{L^2}.
 $$

\subsection {End of the proof}
We have now
\begin{equation}
\label{eq4}
\begin{array}{ll}
\displaystyle{{k\over 8}\int_0^s \sum_{\nu=1}^{+\infty} 2^{2\nu}\|v_\nu(t,\cdot)\|^2_{L^2}\, dt+{\gamma\over 8}\int_0^s \sum_{\nu=0}^{+\infty} \|v_\nu(t,\cdot)\|^2_{L^2}\, dt}\\[0.3 cm]\quad
\displaystyle{\leq  {2\over k}(s+\tau)\sum_{\nu=0}^{+\infty}2^{2\nu}\|v_\nu(s,\cdot)\|^2_{L^2}
+{\gamma \over 2}(s+\tau)\sum_{\nu=0}^{+\infty}\|v_\nu(s,\cdot)\|^2_{L^2}}
\\[0.3 cm]\quad\quad
\displaystyle{+{1\over 2}\tau\,\Phi_\lambda'({\tau\over\beta})\sum_{\nu=0}^{+\infty}\|v_\nu(0,\cdot)\|^2_{L^2}}
\\[0.3 cm]\quad\quad\ 
\displaystyle{-{\gamma\over 8}\int_0^s \sum_{\nu=0}^{+\infty}\|v_\nu(t,\cdot)\|^2_{L^2}\, dt -{5\over 8}\int_0^s (t+\tau) \sum_{\nu=0}^{+\infty}\|\partial_t v_\nu(t,\cdot)\|^2_{L^2}\, dt}
\\[0.3 cm]\quad\quad\ \ 
\displaystyle{-{\alpha k (\log 2)\over 16} \int_0^s(t+\tau)\sum_{\nu=0}^{+\infty}\nu 2^{2\nu} \|v_\nu(t,\cdot)\|^2_{L^2}\, dt}
\\[0.3 cm]\quad\quad\ \ \ 
\displaystyle{+\alpha(\log 2) \tilde Q_1 \int_0^s(t+\tau) \sum_{\nu=0}^{+\infty} 2^{2\nu} \|v_\nu(t,\cdot)\|^2_{L^2}\, dt}
\\[0.3 cm]\quad\quad\ \ \ \ 
\displaystyle{+\int_0^s (t+\tau)(Q_\eta \sum_{\nu=0}^{+\infty}2^{2\nu}\|v_\nu(t,\cdot)\|^2_{L^2}+
\eta \sum_{\nu=0}^{+\infty}\|\partial_t v_\nu(t,\cdot)\|^2_{L^2})\, dt.}
\end{array}
\end{equation}

We choose $\eta$ in such a way that $\eta<{5\over 8}$  and then
$$
-{5\over 8}\int_0^s (t+\tau) \sum_{\nu=0}^{+\infty}\|\partial_t v_\nu(t,\cdot)\|^2_{L^2}\, dt
+\int_0^s (t+\tau) \eta \sum_{\nu=0}^{+\infty}\|\partial_t v_\nu(t,\cdot)\|^2_{L^2}\, dt
\leq 0.
$$

Now if $\nu$ is such that
$$
{\alpha k \log 2\over16}\nu\geq Q_\eta +\alpha(\log 2)\, \tilde Q_1,
$$
then
$$
{\alpha k \log 2\over16}\nu 2^{2\nu} \geq Q_\eta 2^{2\nu} +\alpha(\log 2)\, \tilde Q_1
\,2^{2\nu}.
$$
Consequently, setting $\bar\nu_2:={(16/(\alpha k \log 2))} (Q_\eta +\alpha(\log 2)\, \tilde Q_1)$, we have
$$
\begin{array}{ll}
\displaystyle{-{\alpha k (\log 2)\over 16} \int_0^s(t+\tau)\sum_{\nu= \bar\nu_2}^{+\infty}\nu 2^{2\nu} \|v_\nu(t,\cdot)\|^2_{L^2}\, dt}
\\[0.3 cm]
\quad\qquad\displaystyle{+\alpha(\log 2) \tilde Q_1 \int_0^s(t+\tau) \sum_{\nu=\bar\nu_2}^{+\infty} 2^{2\nu} \|v_\nu(t,\cdot)\|^2_{L^2}\, dt}
\\[0.3 cm]\qquad\qquad\qquad\displaystyle{+Q_\eta \int_0^s (t+\tau)\sum_{\nu=\bar\nu_2}^{+\infty}2^{2\nu}\|v_\nu(t,\cdot)\|^2_{L^2}\, dt\leq 0.}
\end{array}
$$
Finally, eventually choosing a larger $\gamma$  in such a way that
$$
{\gamma\over 8}\geq (\sigma+\tau)(\tilde Q_1\, \alpha(\log 2) +Q_\eta) 2^{2\bar\nu_2},
$$
we obtain
$$
\begin{array}{ll}
\displaystyle{-{\gamma\over 8} \int_0^s \sum_{\nu= 0}^{\bar\nu_2-1} \|v_\nu(t,\cdot)\|^2_{L^2}\, dt +\alpha(\log 2) \tilde Q_1 \int_0^s(t+\tau) \sum_{\nu= 0}^{\bar\nu_2-1} 2^{2\nu} \|v_\nu(t,\cdot)\|^2_{L^2}\, dt}
\\[0.3 cm]\qquad\qquad\qquad\qquad\qquad\qquad\qquad\displaystyle{+Q_\eta \int_0^s (t+\tau)\sum_{\nu= 0}^{\bar\nu_2-1}2^{2\nu}\|v_\nu(t,\cdot)\|^2_{L^2}\, dt\leq 0.}
\end{array}
$$

The inequality (\ref{eq4}) becomes
$$
\begin{array}{ll}
\displaystyle{{k\over 8}\int_0^s \sum_{\nu=1}^{+\infty} 2^{2\nu}\|v_\nu(t,\cdot)\|^2_{L^2}\, dt+{\gamma\over 8}\int_0^s \sum_{\nu=0}^{+\infty} \|v_\nu(t,\cdot)\|^2_{L^2}\, dt}\\[0.3 cm]\quad
\displaystyle{\leq  {2\over k}(s+\tau)\sum_{\nu=0}^{+\infty}2^{2\nu}\|v_\nu(s,\cdot)\|^2_{L^2}
+{\gamma \over 2}(s+\tau)\sum_{\nu=0}^{+\infty}\|v_\nu(s,\cdot)\|^2_{L^2}}
\\[0.3 cm]\qquad\qquad\qquad\qquad\qquad\qquad\qquad\qquad
\displaystyle{+{1\over 2}\tau\,\Phi_\lambda'({\tau\over\beta})\sum_{\nu=0}^{+\infty}\|v_\nu(0,\cdot)\|^2_{L^2}}.
\end{array}
$$
From this,  going back to the function $u_\nu$, we have
$$
\begin{array}{ll}
\displaystyle{{k\over 8}\int_0^s e^{2\gamma t}e^{-2\beta \Phi_\lambda({t+\tau\over \beta})}
\sum_{\nu=1}^{+\infty} 2^{2\nu} 2^{-2\alpha \nu t}\|u_\nu(t,\cdot)\|^2_{L^2}\, dt}\\[0.3 cm]\ 
\displaystyle{+{\gamma\over 8}\int_0^s e^{2\gamma t}e^{-2\beta \Phi_\lambda({t+\tau\over \beta})} \sum_{\nu=0}^{+\infty} 2^{-2\alpha \nu t} \|u_\nu(t,\cdot)\|^2_{L^2}\, dt}
\\[0.3 cm]\quad
\displaystyle{\leq  {2\over k}(s+\tau) e^{2\gamma s}e^{-2\beta \Phi_\lambda({s+\tau\over \beta})}\sum_{\nu=0}^{+\infty}2^{2\nu}2^{-2\alpha \nu s}\|u_\nu(s,\cdot)\|^2_{L^2}}
\\[0.3 cm]\qquad \quad
\displaystyle{
+{\gamma \over 2}(s+\tau)e^{2\gamma s}e^{-2\beta \Phi_\lambda({s+\tau\over \beta})}
\sum_{\nu=0}^{+\infty}2^{-2\alpha \nu s} \|u_\nu(s,\cdot)\|^2_{L^2}}
\\[0.3 cm] \qquad\qquad
\displaystyle{+{1\over 2}\tau\,\Phi_\lambda'({\tau\over\beta})e^{-2\beta \Phi_\lambda({\tau\over \beta})}\sum_{\nu=0}^{+\infty}\|u_\nu(0,\cdot)\|^2_{L^2}},
\end{array}
$$
and (\ref{2}) follows immediately from (\ref{d2}), concluding the proof of Theorem \ref{theorem1}.

\section{Proofs of  Theorems \ref{theorem2} and \ref{theorem3}}

We start with a lemma that will be used in the proof of Theorem 2.
\begin{Lemma}
\label{lemma2}
Let $u\in {\mathcal E}$ be a solution of equation (\ref {1}). Then there exists $\gamma_0>0$ such that if $\gamma>\gamma_0$ then the function 
$E(t):= e^{2\gamma t}\|u(t,\cdot)\|^2_{L^2}$ is (weakly) increasing.
\end{Lemma}

\noindent {\it Proof. } It is sufficient to compute the derivative of $E(t)$. We obtain
$$
\begin{array}{ll}
&\displaystyle{{d\over dt}(e^{2\gamma t}\|u(t,\cdot)\|^2_{L^2})=
2\gamma e^{2\gamma t}\|u(t,\cdot)\|^2_{L^2} + 2e^{2\gamma t}\langle u(t,\cdot), \partial_t u(t,\cdot)\rangle_{L^2}}\\[0.3 cm]
&\qquad=\displaystyle{2\gamma e^{2\gamma t}\|u(t,\cdot)\|^2_{L^2}
+2e^{2\gamma t}\langle A(t,\cdot)\nabla_xu(t,\cdot), \nabla_xu(t,\cdot)\rangle_{(L^2)^n}}
\\[0.3 cm]
&\qquad\qquad\displaystyle{ -2e^{2\gamma t}(\langle B(t,\cdot)\nabla_x u(t,\cdot), u(t,\cdot)\rangle_{L^2}+\langle c(t,\cdot) u(t,\cdot), u(t,\cdot)\rangle_{L^2})}
\\[0.3 cm]
&\qquad\geq \displaystyle{2\gamma e^{2\gamma t}\|u(t,\cdot)\|^2_{L^2}
+2e^{2\gamma t}k \|\nabla_x u(t,\cdot )\|^2_{(L^2)^n}}
\\[0.3 cm]
&\qquad\qquad\displaystyle{-2 e^{2\gamma t}n^{1/2} B \|\nabla_x u(t,\cdot )\|_{(L^2)^n}\|u(t,\cdot)\|_{L^2}-2 e^{2\gamma t}C\|u(t,\cdot)\|^2_{L^2} ,}
\end{array}
$$
and the conclusion follows easily.\qed
\vskip 0.3 cm
\noindent
Let us come to the proof of Theorem \ref{theorem2}. 
Let $\sigma$, $\bar \lambda$, $\alpha$, $\bar \gamma$, $M$ as in Theorem \ref{theorem1}.  We choose $\lambda\geq \bar \lambda$ and $\gamma\geq\max\{\bar \gamma, \gamma_0\}$ where $\gamma_0$ is given by Lemma \ref{lemma2}. We set
$\tau={\sigma\over 2}- 2\bar \sigma$ (we recall that $\bar \sigma\in \,]0,\sigma/4[\,$ and then ${\sigma\over 2}- 2\bar \sigma\in \,]0,\sigma/2[$ ). Then (\ref{2}) gives
$$
\begin{array}{ll}
\displaystyle{\int_0^\sigma e^{2\gamma t}e^{-2\beta \Phi_\lambda({t+\tau\over \beta})}
\|u(t,\cdot)\|^2_{H^{1-\alpha t}}\, dt}
\\[0.3 cm]\quad
\displaystyle{\leq M((\sigma+\tau) e^{2\gamma \sigma}e^{-2\beta \Phi_\lambda({\sigma+\tau\over \beta})}\|u(\sigma,\cdot)\|^2_{L^2}
+\tau\,\Phi_\lambda'({\tau\over\beta})e^{-2\beta \Phi_\lambda({\tau\over \beta})}\|u(0,\cdot)\|^2_{L^2}}),
\end{array}
$$
for all $\beta\geq \sigma+\tau$.
Let  now $s\in [0, \bar\sigma]$. Then $2s+\tau\leq 2\bar \sigma+\tau\leq{ \sigma\over2}<\sigma$ and consequently
$$
\begin{array}{ll}
\displaystyle{\int_s^{2s+\tau} e^{2\gamma t}e^{-2\beta \Phi_\lambda({t+\tau\over \beta})}
\|u(t,\cdot)\|^2_{L^2}\, dt}
\\[0.3 cm]\quad
\displaystyle{\leq M((\sigma+\tau) e^{2\gamma \sigma}e^{-2\beta \Phi_\lambda({\sigma+\tau\over \beta})}\|u(\sigma,\cdot)\|^2_{L^2}
+\tau\,\Phi_\lambda'({\tau\over\beta})e^{-2\beta \Phi_\lambda({\tau\over \beta})}\|u(0,\cdot)\|^2_{L^2}}),
\end{array}
$$
where we have used the fact that $\|u(t,\cdot)\|_{L^2}\leq \|u(t,\cdot)\|_{H^{1-\alpha t}}$.
From Lemma \ref{lemma2} $\|u(t,\cdot)\|_{L^2}$ is increasing. Also the function 
$\Phi_\lambda$ is increasing and consequently the function $t\mapsto 
e^{-2\beta \Phi_\lambda({(t+\tau)/ \beta})}$ is decreasing. We deduce that
$$
\begin{array}{ll}
\displaystyle{ e^{2\gamma s}e^{-2\beta \Phi_\lambda({2s+2\tau\over \beta})}
(s+\tau)\|u(s,\cdot)\|^2_{L^2}}
\\[0.3 cm]\quad
\displaystyle{\leq M((\sigma+\tau) e^{2\gamma \sigma}e^{-2\beta \Phi_\lambda({\sigma+\tau\over \beta})}\|u(\sigma,\cdot)\|^2_{L^2}
+\tau\,\Phi_\lambda'({\tau\over\beta})e^{-2\beta \Phi_\lambda({\tau\over \beta})}\|u(0,\cdot)\|^2_{L^2})}.
\end{array}
$$
Then
$$
\begin{array}{ll}
\displaystyle{ \|u(s,\cdot)\|^2_{L^2}\leq  M ({\sigma+\tau\over \tau}) e^{2\gamma \sigma}\Phi_\lambda'({\tau\over\beta})\Big(e^{2\beta( \Phi_\lambda({\sigma/2+\tau\over \beta})-\Phi_\lambda({\sigma+\tau\over \beta}))}\|u(\sigma,\cdot)\|^2_{L^2}}
\\[0.3 cm]\qquad\qquad\qquad\qquad\qquad\qquad\qquad\qquad\quad
\displaystyle{
+e^{2\beta( \Phi_\lambda({\sigma/2+\tau\over \beta})-\Phi_\lambda({\tau\over \beta}))}\|u(0,\cdot)\|^2_{L^2}}\Big)
\\[0.3 cm]
\displaystyle{\phantom{ \|u(s,\cdot)\|^2_{L^2}}\leq  \tilde M \Phi_\lambda'({\tau\over\beta})e^{2\beta( \Phi_\lambda({\sigma/2+\tau\over \beta})-\Phi_\lambda({\sigma+\tau\over \beta}))}\Big(\|u(\sigma,\cdot)\|^2_{L^2}}
\\[0.3 cm]\qquad\qquad\qquad\qquad\qquad\qquad\qquad\qquad\qquad\qquad\qquad
\displaystyle{
+e^{-2\beta \Phi_\lambda({\tau\over \beta})}\|u(0,\cdot)\|^2_{L^2}\Big),}
\end{array}
$$
where $\tilde M$ depends on $\sigma$, $\tau$, $\gamma$ and $M$. We recall that the function $\Phi_\lambda$ is concave, so that
$$
\Phi_\lambda({\sigma/2+\tau\over \beta})-\Phi_\lambda({\sigma+\tau\over \beta})
\leq \Phi'_\lambda({\sigma+\tau\over \beta})({\sigma/2+\tau\over \beta}-{\sigma+\tau\over \beta})= -\Phi'_\lambda({\sigma+\tau\over \beta}){\sigma\over 2\beta},
$$
and then
$$
\|u(s,\cdot)\|^2_{L^2}\leq  \tilde M \Phi_\lambda'({\tau\over\beta})
e^{-\sigma \Phi'_\lambda({\sigma+\tau\over \beta})}(\|u(\sigma,\cdot)\|^2_{L^2}+e^{-2\beta \Phi_\lambda({\tau\over \beta})}\|u(0,\cdot)\|^2_{L^2}).
$$
From (\ref{a}) we have that 
$$
\Phi'_\lambda({\sigma+\tau\over \beta})=	\psi_\lambda({\sigma+\tau\over \tau}{\tau\over \beta})=\exp({({\sigma+\tau\over \tau})^{-\lambda}-1})\Big(\psi_\lambda({\tau\over \beta})\Big)^{({\sigma+\tau\over \tau})^{-\lambda}}.
$$
Then, setting $\tilde \delta:=({(\sigma+\tau)/ \tau})^{-\lambda}$ we obtain that there exists
$\tilde N>0$ such that
$$
\|u(s,\cdot)\|^2_{L^2}\leq  \tilde M \psi_\lambda({\tau\over\beta})
\exp(-\tilde N(\psi_\lambda({\tau\over \beta}))^{\tilde \delta})(\|u(\sigma,\cdot)\|^2_{L^2}+e^{-2\beta \Phi_\lambda({\tau\over \beta})}\|u(0,\cdot)\|^2_{L^2}).
$$
We choose now $\beta$ in such a way that $e^{-\beta \Phi_\lambda({\tau\over \beta})}=
\|u(0,\cdot)\|^{-1}_{L^2}$ i. e.
$$
{\beta\over \tau} \Phi_\lambda({\tau\over \beta})={1\over \tau}\log \|u(0,\cdot)\|_{L^2}.
$$
 We obtain
$ \beta= \tau \Lambda^{-1}({1\over \tau}\log \|u(0,\cdot)\|_{L^2})$ and then there exists
$\bar \rho>0$ such that if  $\|u(0,\cdot)\|_{L^2}\leq \bar \rho$, then $\beta\geq \sigma+\tau$.

Finally we have
$$
\|u(s,\cdot)\|^2_{L^2}\leq  \tilde{\tilde M} 
\exp(-{\tilde N\over 2}[\psi_\lambda({1\over\Lambda^{-1}({1\over \tau}\log \|u(0,\cdot)\|_{L^2}}) ]^{\tilde \delta})(\|u(\sigma,\cdot)\|^2_{L^2}+1),
$$
and, from (\ref{b}),
$$
\|u(s,\cdot)\|^2_{L^2}\leq  \tilde{\tilde M} 
\exp(-\tilde{\tilde N}[{1\over \tau}\log \|u(0,\cdot)\|_{L^2} ]^{\tilde \delta})(\|u(\sigma,\cdot)\|^2_{L^2}+1).
$$
The inequality (\ref{3}) easy follows, concluding the proof of Theorem \ref{theorem2}. 

To prove Theorem \ref{theorem3} it is 
 sufficient to iterate a finite number of times the local result of Theorem \ref{theorem2} choosing for instance $\bar \sigma= \sigma/8$.

\renewcommand{\theequation}{A.\arabic{equation}}
\section*{Appendix}

In the construction of the following example we will follow closely \cite {Pl} (see also \cite{DSP}). 
Let $A$, $B$, $C$, $J$
be four $C^\infty$ functions defined in $\mathbb R$ with $0\leq
A(s),\ B(s),\ C(s)\\ \leq 1$, $-2\leq J(s)\leq 2$ for all $s\in
{\mathbb R}$ and 
$$
\begin{array}{ll}
\displaystyle{A(s)=1\quad{\rm for\ all}\ s\leq {1\over 5}, }&
\quad\displaystyle{A(s)=0\quad{\rm for\ all}\
s\geq {1\over 4},}\\[0.3 cm]
\displaystyle{B(s)=0\quad{\rm for\ all}\ s\leq 0\ {\rm or}\ s\geq 1, }&\quad
\displaystyle{B(s)=1\quad{\rm for\ all}\ {1\over 6}\leq s\leq {1\over 2},}\\[0.3 cm]
\displaystyle{C(s)=0\quad{\rm for\ all}\ s\leq {1\over 4}, }&
\quad\displaystyle{C(s)=1\quad{\rm for\ all}\
s\geq {1\over 3},}\\[0.3 cm]
\displaystyle{J(s)=-2\quad{\rm for\ all}\ s\leq {1\over 6}\ {\rm or}\ s\geq
{1\over 2}, }&\quad
\displaystyle{J(s)=2\quad{\rm for\ all}\ {1\over 5}\leq s\leq {1\over 3}.}\\
\end{array}
$$ 
Let $(a_n)_n$, $(z_n)_n$ be two real sequences such that
\begin{eqnarray}
\displaystyle{0\leq a_n<a_{n+1}\quad {\rm for\ all}\ n\geq
1}& {\rm and}&\displaystyle{\lim_n a_n=+\infty,}\label
{e1}\\[0.3 cm]
\displaystyle{1\leq z_n<z_{n+1}\quad {\rm for\ all}\ n\geq
1}& {\rm and}&\displaystyle{\lim_n
z_n=+\infty.}\label {e2}
\end{eqnarray} 
Let us define $r_n=a_{n+1}-a_n$, $q_1=0$,
$q_n=\sum_{k=2}^nz_k r_{k-1}$ for all
$n\geq 2$, and $p_n=(z_{n+1}-z_n)r_n$. We suppose that
\begin{eqnarray}
\displaystyle{r_n<1\quad {\rm for\ all}\ n\geq
1},\label
{e3}\\[0.3 cm]
\displaystyle{p_n>1\quad {\rm for \ all}\ n\geq 1}.\label {e4}
\end{eqnarray}
We set $A_n(t)=A({t-a_n\over r_n})$, $B_n(t)=B({t-a_n\over
r_n})$,
$C_n(t)=C({t-a_n\over r_n})$ and $J_n(t)=J({t-a_n\over r_n})$. We define
$$
\begin{array}{l}
\displaystyle{v_n(t,x_1)=\exp (-q_n-z_n(t-a_n))\cos \sqrt {z_n} x_1,}\\[0.3 cm]
\displaystyle{w_n(t,x_2)=\exp (-q_n-z_n(t-a_n)+J_n(t)p_n)\cos \sqrt {z_n}
x_2,}\\
\end{array}
$$ 
and, for $n_0\geq 1 $ to be chosen, 
$$ 
u(t,x_1, x_2)= v_{n_0}(t, x_1)
$$ 
for all $t\leq a_{n_0}$, $(x_1, x_2)\in {\mathbb R}^2$ and 
$$ 
u(t,x_1, x_2)= A_n(t)v_n(t,x_1)+B_n(t)w_n(t,x_2)+C_n(t)v_{n+1}(t,x_1)
$$
for all $n\geq n_0$,  $a_n\leq t\leq a_{n+1}$ and $(x_1, x_2)\in {\mathbb R}^2$.
If, for all $\alpha$, $\beta$ $\gamma>0$,
\begin{equation}
\lim_n\exp(-q_n+2p_n)z_{n+1}^\alpha p_n^\beta r_n^{-\gamma}=0,
\label{e5}
\end{equation} 
then $u$ is a $C^\infty_B({\mathbb R}^3)$ function, where $C^\infty_B$ denotes the smooth functions which are bounded with bounded derivatives. We define
$$ l(t)=\left\{
\begin{array}{ll}
\displaystyle{1} &\displaystyle{\quad{\rm for \ all}\ t\leq a_1,
}\\[0.3 cm]
\displaystyle{1+J'_n(t)p_nz^{-1}_n}&\displaystyle{\quad{\rm for\ all }\ a_n\leq t\leq
a_{n+1} .}\\
\end{array}
\right.
$$ The condition
\begin{equation}
\sup_{n\geq n_0}\; \{p_nr_n^{-1}z_n^{-1}\}\leq {1\over 2\|J'\|_{L^\infty}}
\label{e6}
\end{equation} 
guarantees that the operator $\text{\it
\L}=\partial_t-\partial^2_{x_1}-l(t)\partial^2_{x_2}$ is parabolic.
The
 function $l$ is smooth and it is Log-Lipschitz continuous on $\mathbb R$ (i. e.   
 $$
 \sup_{t_1<t_2}{|l(t_2)-l(t_1)|\over |t_2-t_1|(|\log|t_2-t_1||+1)}<+\infty)
 $$
under the following condition
\begin{equation}
\sup_{n}\; \{{p_nr_n^{-1}z_n^{-1}\over r_n\log ({1\over r_n})}\}<+\infty.
\label{e7}
\end{equation} 
Finally we define
$$
\begin{array}{l}
\displaystyle{b_1=- {\text{\it \L} u\over
u^2+(\partial_{x_1}u)^2+(\partial_{x_2}u)^2}\partial_{x_1}u,}\\[0.3cm]
\displaystyle{b_2=- {\text{\it \L} u\over
u^2+(\partial_{x_1}u)^2+(\partial_{x_2}u)^2}\partial_{x_2}u,}\\[0.3cm]
\displaystyle{c=- {\text{\it \L} u\over
u^2+(\partial_{x_1}u)^2+(\partial_{x_2}u)^2}u}\\
\end{array}
$$ 
and, as in \cite {Pl}, the coefficients $b_1$, $b_2$, $c$ will be
in $C^\infty_B({\mathbb R}^3)$ if, for all $\alpha$, $\beta$, $\gamma>0$,
\begin{equation}
\lim_n\exp(-p_n)z_{n+1}^\alpha p_n^\beta r_n^{-\gamma}=0.
\label{e8}
\end{equation} 
We set
$$
a_1=0,\ a_n= \sum _{j=2}^{n}{1\over j\log j}\quad{\rm for\ all}\ n\geq 2,
$$
and 
$$
z_n= n^4 \quad{\rm for\ all}\ n\geq 1.
$$
With these choices the conditions (\ref {e1}), (\ref {e2}), (\ref {e3}) and (\ref {e4}) are trivial and also  (\ref {e5}), (\ref {e7}) and (\ref {e8}) are easily verified. From (\ref {e7}) and the fact that $\lim_{n}r_n=0$ we deduce that 
$$
\lim_{n}p_nr_n^{-1}z_n^{-1}=0,
$$
consequently it is possible to choose $n_0$ such that (\ref {e6}) is verified.

Let now 
$$
n_{1,k}=\big[\exp(\exp( k))\big]+2\quad{\rm and}\quad n_{2,k}=\big[\exp(\exp (k+{1\over k})\big]+1,
$$
where $\big[x\big]$ denotes the integer part of $x$ and where $k$ is taken in such a way that $n_{1,k}\geq n_0$. We fix, for $h=1,2$,
$$
t_{h,k}=a_{n_{h,k}}=\sum_{j=1}^{n_{h,k}-1}r_j, 
$$
We have $\lim_k t_{1,k}=+\infty$ and 
$$
t_{2,k}- t_{1,k}=\sum_{j=n_{1,k}}^{n_{2,k}-1}r_j\leq 
\int_{n_{1,k}-1}^{n_{2,k}-1}{1\over x\log x}dx\leq 
\int_{\exp(\exp( k))}^{\exp(\exp (k+{1\over k})}{1\over x\log x}dx={1\over k},
$$
so that $\lim_k t_{2,k}- t_{1,k}=0$. Our intent is to prove that 
\begin{equation}
\lim_k {\exp(-q_{n_{1,k}})\over  \exp(-\delta q_{n_{2,k}})}=+\infty
\label{e9}
\end{equation}
for all $\delta\in (0,1)$. 
Since, by the choice of $(a_n)_n$ and $(z_n)_n$ we have that $q_n=\sum_{j=2}^n
{j^3 /\log j}$, it is immediate to obtain that
$$
q_{n_{2,k}}\geq q_{n_{1,k}}+{n_{1,k}^3\over \log (n_{1,k})}(n_{2,k}-n_{1,k}),
$$
and (\ref{e9}) is a consequence of
$$
\lim_k \delta(q_{n_{1,k}}+{n_{1,k}^3\over \log (n_{1,k})}(n_{2,k}-n_{1,k}))-q_{n_{1,k}}=+\infty.
$$
This result will be implied by 
\begin{equation}
\lim_k \delta{n_{1,k}^3\over \log (n_{1,k})}(n_{2,k}-n_{1,k})-q_{n_{1,k}}=+\infty,
\label{e10}
\end{equation}
 for all $\delta\in (0,1)$. 
Easily we have that $q_{n_{1,k}}\leq n_{1,k}^4$ for all $k$ and consequently (\ref{e10})
may be deduced from 
$$
\lim_k \delta'{n_{1,k}^3 n_{2,k}\over \log (n_{1,k})}-n_{1,k}^4=+\infty
$$
i. e. 
$$
\lim_k \delta'{n_{2,k}\over \log (n_{1,k})}-n_{1,k}=+\infty\quad {\rm for\ all} \ \delta'\in (0,1),
$$
which can be elementary obtained substituting $n_{1,k}$ and $n_{2,k}$ with their values.
Summing up we have the following result.

\begin{Theorem}
\label{te}
There exist
\begin{itemize}
\item $l\in C^\infty({\mathbb R})$, $l$ Log-Lipschitz continuous, 
 $1/2\leq l(t)\leq 3/2$ for all $t\in {\mathbb R}$,
\item  $b_1$, $b_2$, $c$ and $u\in C_B^\infty({\mathbb R^3})$, $2\pi$-periodic  with respect to $x_1$ and $x_2$, 
\item $(t_{1,n})_n$, $(t_{2,n})_n$ increasing sequences in $\mathbb R$, $1>t_{2,n}-t_{1,n}>0$ for all $n$, $\lim_n t_{1,n}=+\infty$ and $\lim_n t_{2,n}-t_{1,n}=0$,
\end{itemize}
such that
$$
\partial_tu-\partial^2_{x_1}u-l\partial^2_{x_2}u+b_1\partial_{x_1}u+b_2\partial_{x_2}u
+cu= 0
$$
for all $(t,x_1, x_2)\in {\mathbb R}^3$ and 
$$
\lim_n {\|u(t_{1, n}, \cdot, \cdot)\|_{L^2([0,2\pi] \times [0,2\pi])}\over
\|u(t_{2, n}, \cdot, \cdot)\|^\delta_{L^2([0,2\pi] \times [0,2\pi])}}=+\infty
$$
for all $\delta\in (0,1)$.
\end{Theorem}

We define now, for $(t,x_1,x_2)\in [0,1]\times {\mathbb R}^2$, 
$$
\begin{array}{l}
 \displaystyle{l_n(t)=l(t_{2,n}-t)},\\[0.3 cm]
 \displaystyle{u_n(t,x_1,x_2)= u(t_{2,n}-t, x_1, x_2)},
\end{array}
$$
and similarly for $b_{1,n}$, $b_{2,n}$ $c_n$. We set $t_n=t_{2,n}-t_{1,n}$ and 
$$
L_n= \partial_t +\partial^2_{x_1}+l_n\partial^2_{x_2}u-b_{1, n}\partial_{x_1}-b_{2,n }\partial_{x_2}
-c_n.
$$
We have that $(L_n)_n$ is a sequence of uniformly backward parabolic operators with uniformly Log-Lipschitz continuous coefficients in the principal part and uniformly bounded coefficients in lower order terms. $(u_n)_n$ is a sequence of smooth uniformly bounded solutions of $L_nu_n=0$ on $[0,1]\times {\mathbb R}^2$, with 
$$
\lim_n \|u_n(0,  \cdot,\cdot)\|_{L^2([0,2\pi]\times [0,2\pi])}=0.
$$ 

We have that $\lim_n t_n=0$ and
$$
\lim_n {\|u_n(t_n, \cdot, \cdot)\|_{L^2([0,2\pi] \times [0,2\pi])}\over
\|u_n(0, \cdot, \cdot)\|^\delta_{L^2([0,2\pi] \times [0,2\pi])}}=+\infty
$$
for all $\delta\in (0,1)$: {\it it is not possible to obtain a result similar to that of Hurd \cite{H} or Agmon and Nirenberg \cite{AN} if Lipscihtz continuity is replaced by Log-Lipschitz continuity}.

\end{document}